\documentclass[10pt]{article}

\usepackage[usenames]{color}
\usepackage{amssymb,psfig,epsfig,latexsym,graphicx}
\usepackage{amsmath}

\definecolor{webgreen}{rgb}{0,.5,0}
\definecolor{webbrown}{rgb}{.6,0,0}

\newtheorem{theorem}{Theorem}

\newtheorem{coro}[theorem]{Corollary}

\newcommand{\eqn}[1]{(\ref{#1})}
\newcommand{\bsq}{{\vrule height .9ex width .8ex depth -.1ex }}

\newcommand{\ON}{{\tt ON}}
\newcommand{\OFF}{{\tt OFF}}
\newcommand{\DEAD}{{\tt DEAD}}
\newcommand{\F}{{\mathbf F}}

\newcommand{\al}{\alpha}
\newcommand{\be}{\beta}
\newcommand{\ga}{\gamma}
\newcommand{\de}{\delta}
\newcommand{\eps}{\varepsilon}
\newcommand{\beql}[1]{\begin{equation}\label{#1}}
\newcommand{\eeq}{\end{equation}}
\DeclareMathOperator{\wt}{wt}
\newcommand{\mydot}{\mathord{\cdot}}
\newcommand{\sA}{{\mathcal A}}
\newcommand{\sC}{{\mathcal C}}
\newcommand{\sQ}{{\mathcal Q}}

\newcommand{\sT}{{\mathcal T}}

\newcommand{\RR}{{\mathbb R}}

\newcommand{\ZZ}{{\mathbb Z}}

\makeatletter
\def\@sect#1#2#3#4#5#6[#7]#8{\ifnum #2>\c@secnumdepth
     \def\@svsec{}\else
     \refstepcounter{#1}\edef\@svsec{\csname the#1\endcsname.\hskip .75em }\fi
     \@tempskipa #5\relax
      \ifdim \@tempskipa>\z@
        \begingroup #6\relax
          \@hangfrom{\hskip #3\relax\@svsec}{\interlinepenalty \@M #8\par}%
        \endgroup
       \csname #1mark\endcsname{#7}\addcontentsline
         {toc}{#1}{\ifnum #2>\c@secnumdepth \else
                      \protect\numberline{\csname the#1\endcsname}\fi
                    #7}\else
        \def\@svsechd{#6\hskip #3\@svsec #8\csname #1mark\endcsname
                      {#7}\addcontentsline
                           {toc}{#1}{\ifnum #2>\c@secnumdepth \else
                             \protect\numberline{\csname the#1\endcsname}\fi
                       #7}}\fi
     \@xsect{#5}}
\def\@begintheorem#1#2{\it \trivlist \item[\hskip \labelsep{\bf #1\ #2.}]}

\def\section{\@startsection {section}{1}{\z@}{-3.5ex plus -1ex minus
 -.2ex}{2.3ex plus .2ex}{\normalsize\bf}}
\makeatother
\makeatletter
\def\subsection{\@startsection {subsection}{1}{\z@}{-3.5ex plus -1ex minus
 -.2ex}{2.3ex plus .2ex}{\normalsize\bf}}

\makeatother

\begin{document}

%%%%%%%%%%%%%%%%%%%%%%%%%%%%%%%%%%%%%%
% Title, authors
%%%%%%%%%%%%%%%%%%%%%%%%%%%%%%%%%%%%%%

\begin{center}
%{\large\bf The Toothpick Sequence } \\
{\large\bf The Toothpick Sequence and Other Sequences from Cellular Automata }\\
\vspace*{+.2in}

David Applegate, \\
AT\&T Shannon Labs, \\
180 Park Ave., Florham Park, \\
NJ 07932-0971, USA, \\
Email: david@research.att.com, \\

\vspace*{+.1in}

Omar E. Pol, \\
Nazca 5482, CP 1419, \\
Buenos Aires, ARGENTINA, \\
Email: info@polprimos.com, \\

\vspace*{+.1in}

N. J. A. Sloane${}^{(a)}$, \\
AT\&T Shannon Labs, \\
180 Park Ave., Florham Park, \\
NJ 07932-0971, USA, \\
Email: njas@research.att.com. \\

\vspace*{+.2in}
${}^{(a)}$ To whom correspondence should be addressed.

\vspace*{+.2in}
February 13, 2010; revised April 21, 2010, October 2, 2010 \\
\vspace*{+.2in}

%%%%%%%%%%%%%%%%%%%%%%%%%%%%%%%%%%%%%%
% Abstract
%%%%%%%%%%%%%%%%%%%%%%%%%%%%%%%%%%%%%%

{\bf Abstract}
\end{center}

A two-dimensional arrangement of toothpicks is constructed
by the following iterative procedure.
At stage $1$, place a single toothpick of length $1$ on a square grid,
aligned with the $y$-axis.
At each subsequent stage, for every exposed toothpick end,
place a perpendicular toothpick centered at that end.
The resulting structure has a fractal-like appearance.
We will analyze the {\em toothpick sequence},
which gives the total number of toothpicks after $n$ steps.
We also study several related sequences that
arise from enumerating active cells in cellular automata.
Some unusual recurrences appear: a typical example
is that instead of the Fibonacci recurrence, which
we may write as $a(2+i) = a(i) + a(i+1)$,
we set $n = 2^k+i$ ($0 \le i < 2^k$), and then
$a(n)=a(2^k+i)=2a(i)+a(i+1)$. 
The corresponding  generating functions look like 
$\prod_{k \ge 0} (1+x^{2^k-1}+2x^{2^k})$ and variations thereof.

\vspace{0.8\baselineskip}
Keywords:  cellular automata (CA), enumeration, Holladay-Ulam CA,
Schrandt-Ulam CA, Ulam-Warburton CA, Rule 942, 
Sierpi\'{n}ski triangle

\vspace{0.8\baselineskip}
AMS 2000 Classification: Primary 11B85

%%%%%%%%%%%%%%%%%%%%%%%%%%%%%%%%%%%%%%
% Section: introduction
%%%%%%%%%%%%%%%%%%%%%%%%%%%%%%%%%%%%%%

\section{Introduction}\label{SecTP}

We start with an infinite sheet of graph paper and an infinite 
supply of line segments of length $1$, called ``toothpicks.''
At stage $1$, we place a toothpick on the $y$-axis
and centered at the origin.
Each toothpick we place has two ends, and an end is said to be ``exposed''
if this point on the plane is neither the end nor the midpoint
of any other toothpick.

At each subsequent stage, for every exposed toothpick end, we place 
a toothpick centered at that end and
perpendicular to that toothpick.
The toothpicks placed at odd-numbered stages are therefore
all parallel to the $y$-axis, while
those placed at even-numbered stages are parallel to the $x$-axis.

\begin{figure}[htb]
%\centerline{\includegraphics[width=5in]{toothFig1.pdf}}
\centerline{\includegraphics[width=4.5in]{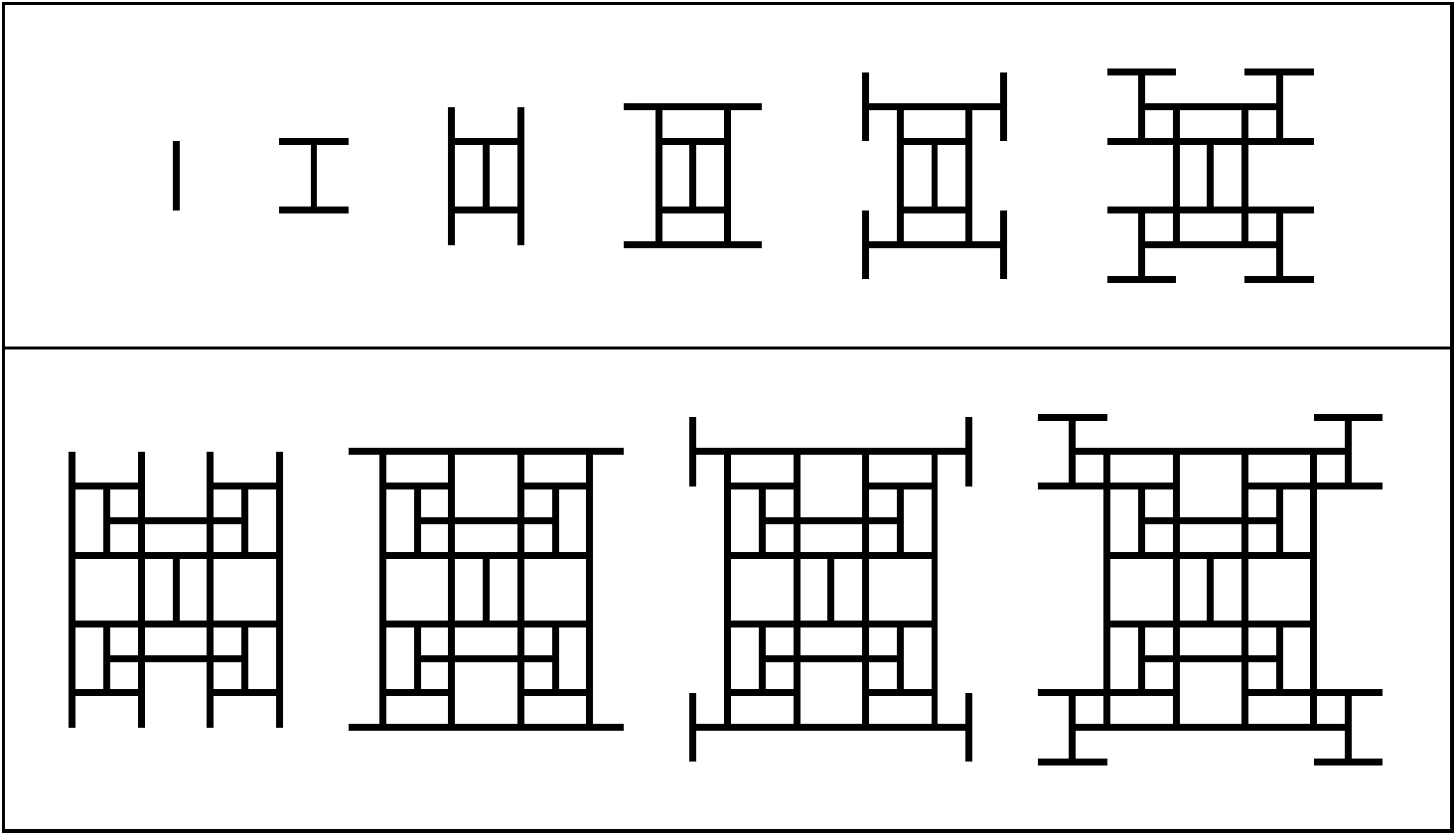}}
\caption{First ten stages of the 
evolution of the toothpick structure.
The numbers of toothpicks in the successive stages, $T(1), \ldots,
T(10)$, are $1, 3, 7, 11, 15, 23, 35, 43, 47, 55$.}
\label{FigTP1}
\end{figure}

\begin{figure}[htb]
%\centerline{\includegraphics[width=5in]{A139250_53_500_0_0_1.pdf}}
\centerline{\includegraphics[width=3.59in]{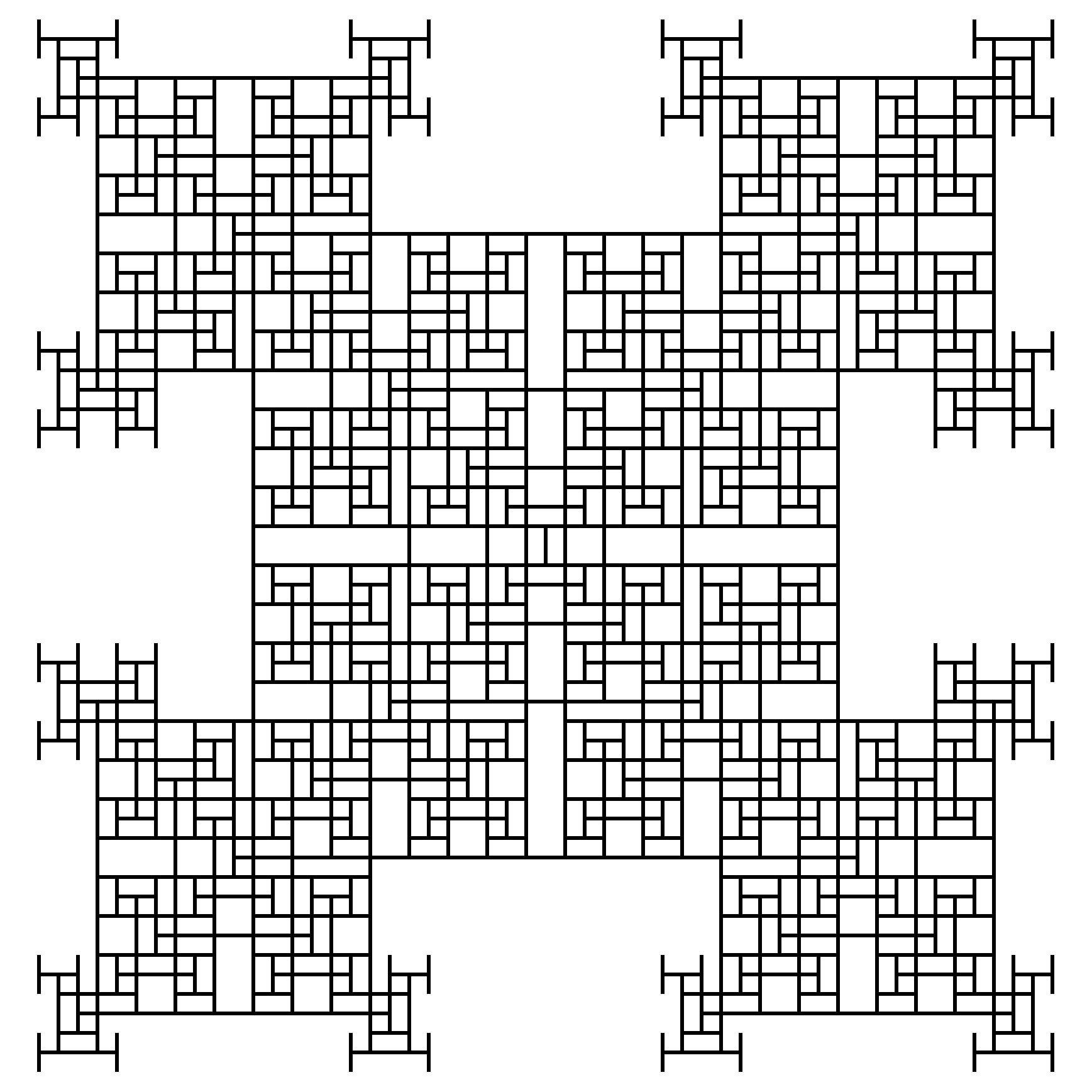}}
\caption{ The toothpick structure after 53 stages 
(there are $T(53)=1379$ toothpicks).}
\label{Fig53}
\end{figure}

\begin{figure}[htb]
%\centerline{\includegraphics[width=5in]{A139250_64_500_0_0_1.pdf}}
\centerline{\includegraphics[width=4.4in]{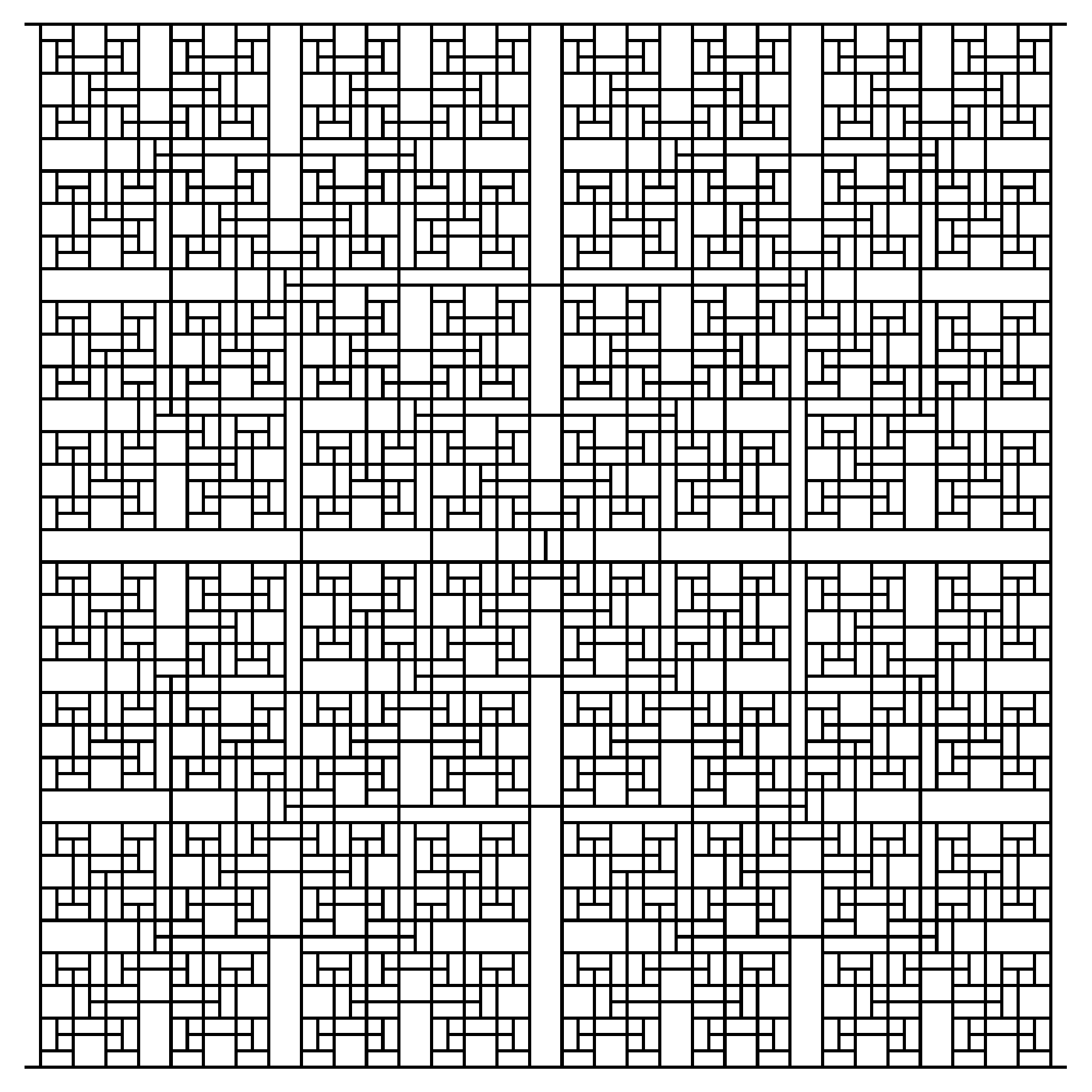}}
\caption{ The toothpick structure after 64 stages
(there are $T(64)=2731$ toothpicks).}
\label{Fig64}
\end{figure}

Fig. \ref{FigTP1} shows the first ten stages of the
evolution of the toothpick structure and
Figs. \ref{Fig53}, \ref{Fig64} show the structure
after respectively $53$ and $64$ stages.
%to illustrate the fractal-like behavior.

Let $t(n)$ ($n \ge 1$) denote the number of toothpicks added
at the $n$th stage, with $t(0)=0$, and let 
$T(n) := \sum_{i=0}^{n} t(i)$ be the total number of
toothpicks after $n$ stages.
The initial values of $t(n)$ and $T(n)$ are shown in Table \ref{TabTP1}.
These two sequences respectively
form entries A139251 and A139250 in \cite{OEIS}.

\begin{table}[htb]
$$
\begin{array}{|c|rrrrrrrrrr|} 
\hline
n   & 0 & 1 & 2 & 3 & 4 & 5 & 6 & 7 & 8 & 9  \\
t(n) & 0 & 1 & 2 & 4 & 4 & 4 & 8 & 12 & 8 & 4 \\
T(n) & 0 & 1 & 3 & 7 & 11 & 15 & 23 & 35 & 43 & 47 \\
\hline
n & 10 & 11 & 12 & 13 & 14 & 15  & 16 & 17 & 18 & 19  \\
t(n) & 8 & 12 & 12 & 16 & 28 & 32 & 16 & 4 & 8 & 12 \\
T(n) & 55 & 67 & 79 & 95 & 123 & 155 & 171 & 175 & 183 & 195 \\ 
\hline
n & 20 & 21 & 22 & 23 & 24 & 25 & 26 & 27 & 28 & 29  \\
t(n) & 12 & 16 & 28 & 32 & 20 & 16 & 28 & 36 &  40 & 60  \\
T(n) & 207 & 223 &  251 & 283 & 303 & 319 & 347 & 383 & 423 & 483  \\
\hline
n & 30 & 31 & 32 & 33 & 34 & 35 & 36 & 37 & 38 & 39  \\
t(n) & 88 & 80 & 32 & 4 & 8 & 12 & 12 & 16 & 28 & 32  \\
T(n) & 571 & 651 & 683 & 687 & 695 & 707 & 719 & 735 & 763 &  795  \\
\hline
n & 40 & 41 & 42 & 43 & 44 & 45 & 46 & 47 & 48 & 49 \\
t(n) & 20 & 16 & 28 & 36 & 40 & 60 & 88 & 80 & 36 & 16 \\
T(n) & 815 & 831 & 859 & 895 & 935 & 995 & 1083 & 1163 & 1199 & 1215 \\
\hline
\end{array}
$$
\caption{The toothpick sequences $t(n)$ and $T(n)$ for $0 \le n \le 49$.}
\label{TabTP1}
\end{table}

The first question is, what are
the numbers $t(n)$ and $T(n)$?
We start by finding recurrences that they satisfy.
As the number of stages grows, the array of toothpicks 
has a recursive, fractal-like structure, as suggested by Figs.
\ref{FigTP1}, \ref{Fig53}, \ref{Fig64}.
(For a dramatic illustration of the fractal structure, see the movie 
linked to entry A139251 in \cite{OEIS}.)
In order to analyze this structure, we consider a variant,
the ``corner'' sequence, which starts from a half-toothpick
protruding from one quadrant of the plane.
In \S\ref{SecCS} we establish a recurrence for
the corner sequence (Theorem \ref{th1})
and in \S\ref{SecTP2} we use this to find recurrences for 
$t(n)$ and $T(n)$ (Theorem \ref{th2}, Corollary \ref{corT}).
Section \S\ref{SecWin} gives a similar recurrence for the number
of squares and rectangles that are created in the toothpick
structure at the $n$th stage, and \S\ref{SecFrac}
gives a more precise description of the fractal-like
behavior and discusses the asymptotic growth of $T(n)$.

The recurrences make it easy to compute a large number
of values of $t(n)$ and $T(n)$, so in a sense the initial problem
has now been solved.

However, the toothpick structure is reminiscent of another,
simpler, two-dimensional structure, the
arrangement of square cells produced by the {\em Ulam-Warburton
cellular automaton} (Ulam \cite{Ulam62}, Singmaster \cite{Sing03},
Stanley and Chapman \cite{Stan94}, Wolfram \cite[p.~928]{NKS}).
For this structure there is an explicit formula for
the number of $\ON$ cells at the $n$th
stage, and a simple generating function for these numbers,
as we will see in \S\ref{SecUW}, Theorem \ref{th3}.

This hint led us to look for
a similar generating function 
and an explicit formula for the toothpick sequence.
Our first attempt was a failure, but provided an
surprising connection  with the Sierpi\'{n}ski triangle,
described in \S\ref{SecLT}.

The generating functions for the toothpick sequence and for
a number of related sequences have an interesting form:
they can be written as
\beql{EqGF1}
x(\al + \be x) \prod_{k \ge \eps} (1 + \ga x^{2^k - 1} + \de x^{2^k}) \, ,
\eeq 
for appropriate integers $\al, \be, \ga, \de, \eps$.
In \S\ref{SecGF} we describe the relationship between such
generating functions and recurrences for the underlying sequence
(Theorem \ref{th4}).
The generating function for the toothpick sequence is then 
established in Theorem \ref{th5}.

Generating functions of the form $\prod_{k \ge 1} (1+g_k x^k)$
have been used in combinatorics and number theory
for a long time (for a survey see \cite{GGM88}),
but generating functions of the form \eqn{EqGF1}
may be new---at least, until the commencement of this work,
there were essentially no examples among the 170,000 entries in \cite{OEIS}.

The following section, \S\ref{SecEX},
gives a general method for obtaining explicit formulas 
from the generating functions (Theorem \ref{th6}),
and the particular formula for the toothpick sequence 
is given in Theorem \ref{th7}.

%We discuss some further sequences of this type in 
%Section \ref{SecGFG}. 

Both the toothpick structure and the Ulam-Warburton structure
are examples of cellular automata defined on graphs,
and we discuss this general framework in \S\ref{SecCAG}.
We have not been able to find much earlier work
on the enumeration of active cells in cellular automata---the
Stanley and Chapman {\em American Mathematical Monthly} problem 
\cite{Stan94} and the Singmaster article \cite{Sing03}
being exceptions.
We would appreciate hearing of any references we have overlooked.

Of course, two well-known examples show that one cannot hope
to enumerate the active states in arbitrary cellular automata:
the one-dimensional cellular automaton defined by Wolfram's ``Rule~$30$'' 
\cite{Weiss}, \cite{Wolf83}, \cite{NKS}
behaves chaotically, and the two-dimensional cellular automaton
corresponding to Conway's ``Game of Life''
\cite{WW2001}
is a universal Turing machine.

However, we were able to apply our techniques (with
varying degrees of success) to a number of other 
cellular automata defined on graphs,
and the last four sections discuss some of these.
Section \ref{SecTTP} discusses a structure built using T-shaped toothpicks.
Sections \ref{SecMC}-\ref{Sec8N} discusses variations on the 
Ulam-Warburton cellular automaton.
Section \ref{SecMC} deals with the ``Maltese cross'' 
or Holladay-Ulam structure studied in \cite{Ulam62},
as well as some other structures mentioned in that paper. 
Section \ref{Sec942} considers what happens if we change the rule for the Ulam-Warburton
cellular automaton of \S\ref{SecUW} so that a cell is turned $\ON$
if and only if one or four of its neighbors is $\ON$. 
The final section (\S\ref{Sec8N}) discusses
what happens if we change the definition of the Ulam-Warburton
cellular automaton to allow all eight
neighbors of a square to affect the next stage.
There is another variation that could have been included here,
in which the rule is that a cell {\em changes state}
if exactly one of its four neighbors is $\ON$.
Again we have a formula for the number of $ON$ cells after $n$ 
generations---see entries A079315, A079317 in \cite{OEIS}.
Many further examples of sequences based on generalized
toothpick structures and cellular automata are listed in \cite{toothlist}.

\paragraph{Notation.}
Our cellular automata are synchronous, and we normally use
the symbol $n$ to index the successive stages.
Cells are either $\ON$ or $\OFF$.
In all the examples we consider here, once a cell
is $\ON$ it stays $\ON$.
Lower case letters (e.g. $a(n)$) will denote
the number of toothpicks added, or cells whose state is
changed from $\OFF$ to $\ON$,
at the $n$th stage, and  the corresponding
upper case letters (e.g. $A(n)$) will denote
the total number of toothpicks or $\ON$ cells after
$n$ stages (the partial sums of the $a(n)$).
By the {\em generating function} for a sequence $a(n)$ (say),
we will always mean the ordinary generating 
function $\sA(x) := \sum_{n=0}^{\infty} a(n) x^n$.
If $\sA(x)$ is the generating function for $a(n)$, then
$\sA(x)/(1-x)$ is the generating function for $A(n)$.

\paragraph{Remarks.}

1. A common dilemma in combinatorics is 
whether to index the first counting step with $n=0$
or $n=1$. In this paper we have consistently 
started the enumerations with zero objects (toothpicks or
$\ON$ cells) at stage $0$, adding the initial object at stage $1$.
This seems natural, and is the indexing used for most
of these sequences in \cite{OEIS}.
On the other hand, this is responsible for the leading factor of $x$ in
the generating functions \eqn{EqGF1}, \eqn{EqGF2}, \eqn{EqGCS1}, etc.,
and for the fact that in the recurrences \eqn{EqCS1}, \eqn{EqTP1}, etc.,
the exceptional cases occur at the beginning of each
block of $2^k$ terms, rather than at the end. If they had occurred at
the ends of the blocks, the beginnings of 
all the blocks would have agreed,
which would have made the triangular arrays such as that in Table
\ref{TabTPtri} look rather nicer (compare Table \ref{TabTPtri2}).
Probably there is no perfect solution to the problem,
and so we have followed the indexing used in \cite{OEIS}.

2. Reference \cite{OEIS} contains several hundred sequences related to
the toothpick problem, many more than could be mentioned here. 
For a full list see \cite{toothlist} and also
the entries in the index to \cite{OEIS}
under ``cellular automata.''

3. The computer language {\em Mathematica{\textregistered}} \cite{Wolf96}
has a collection of commands that can often be used to 
display structures produced by cellular automata and to count the $\ON$ states. 
For example, the command

{\tt Map[Function[Apply[Plus,Flatten[\#1]]],CellularAutomaton[\{}\\
{\tt \hspace*{0.3in} 686,\{2,\{\{0,2,0\},\{2,1,2\},\{0,2,0\}\}\},\{1,1\}\},\{\{\{1\}\},0\},200]]}

\noindent{}produces the first 200 terms of the sequence $U(n)$
giving the number of $\ON$ states in the
Ulam-Warburton cellular automaton discussed in \S\ref{SecUW}.

%%%%%%%%%%%%%%%%%%%%%%%%%%%%%%%%%%%%%%
% Section: 
%%%%%%%%%%%%%%%%%%%%%%%%%%%%%%%%%%%%%%

\begin{figure}[htbp]
\centerline{\includegraphics[width=4.5in]{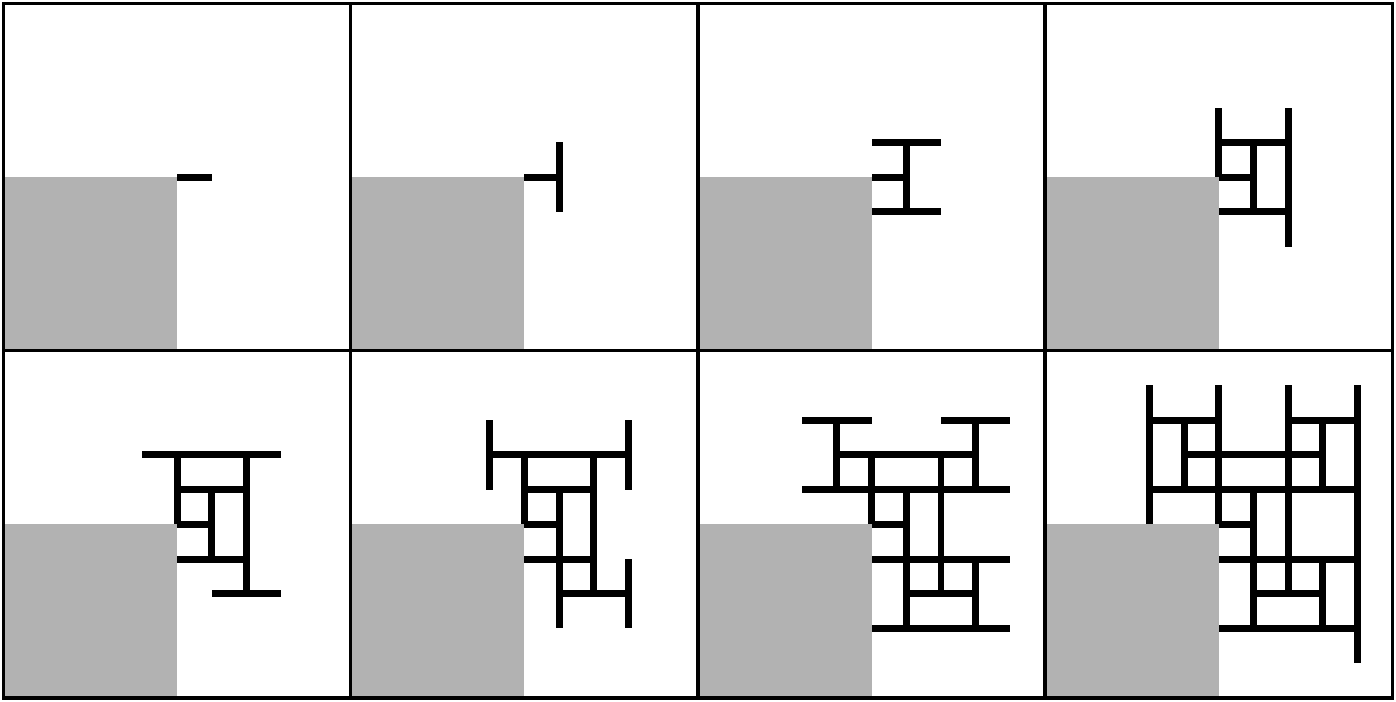}}
\caption{Stages $0$ through $7$ of the
evolution of the corner structure.
The numbers of toothpicks in the successive stages, $C(0), \ldots,
C(7)$, are $0, 1, 3, 6, 9, 13, 20, 28$.}
\label{FigCS1}
\end{figure}

\section{ The corner sequence}\label{SecCS}

In order to understand the toothpick structure, 
it is helpful to first consider what happens if
one quadrant of the plane is excluded.
We impose the rule that no toothpick may cross into the third quadrant
of the plane, and only ends of toothpicks may touch
the negative $x$- or $y$-axes.
At stage $0$, we place a half-toothpick extending horizontally
from the origin to the point $(\frac{1}{2}, 0)$.
The structure is then allowed to grow using the rule 
for the original toothpick sequence.
The corner sequence is relevant because it describes
how the main toothpick structure grows.

Let $c(n)$ ($n\ge 1$) denote the number of
toothpicks added at the $n$th stage, with $c(0)=0$,
and let $C(n) := \sum_{i = 0}^{n} c(i)$ be the total number
of toothpicks after $n$ stages.
These are respectively entries A152980 and A153006 in \cite{OEIS}.

Fig. \ref{FigCS1} shows stages $0$ through $7$ of the
evolution of the corner structure.
Note that the first toothpick added, at stage $1$
(with midpoint at the end of the initial half-toothpick),
matches the initial toothpick of the original toothpick 
sequence, except that it is shifted a half-unit to the left.
The initial values of $c(n)$ and $C(n)$ are shown in Table \ref{TabCS1}.

\begin{table}[htbp]
$$
\begin{array}{|c|rrrrrrrrrr|} 
\hline
n   & 0 & 1 & 2 & 3 & 4 & 5 & 6 & 7 & 8 & 9  \\
c(n) & 0 & 1 & 2 & 3 & 3 & 4 & 7 & 8 & 5 & 4 \\
C(n) & 0 & 1 & 3 & 6 & 9 & 13 & 20 & 28 & 33 & 37 \\
\hline
n & 10 & 11 & 12 & 13 & 14 & 15 & 16 & 17 & 18 & 19 \\
c(n) & 7 & 9 & 10 & 15 & 22 & 20 & 9 & 4 & 7 & 9 \\
C(n) & 44 & 53 & 63 & 78 & 100 & 120 & 129 & 133 & 140 & 149 \\
\hline
n & 20 & 21 & 22 & 23 & 24 & 25 & 26 & 27 & 28 & 29 \\
c(n) & 10 & 15 & 22 & 21 & 14 & 15 & 23 & 28 & 35 & 52 \\
C(n) & 159 & 174 & 196 &  217 & 231 & 246 & 269 & 297 & 332 & 384 \\
\hline
n & 30 & 31 & 32 & 33 & 34 & 35 & 36 & 37 & 38 & 39 \\
c(n) & 64 & 48 & 17 & 4 & 7 & 9 & 10 & 15 & 22 & 21 \\
C(n) & 448 & 496 & 513 & 517 & 524 & 533 & 543 & 558 & 580 & 601 \\
\hline
\end{array}
$$
\caption{The corner sequences $c(n)$ and $C(n)$ for $0 \le n \le 39$.}
\label{TabCS1}
\end{table}

\begin{figure}[htbp]
%\centerline{\includegraphics[width=5in]{A153006.pdf}}
\centerline{\includegraphics[width=4.5in]{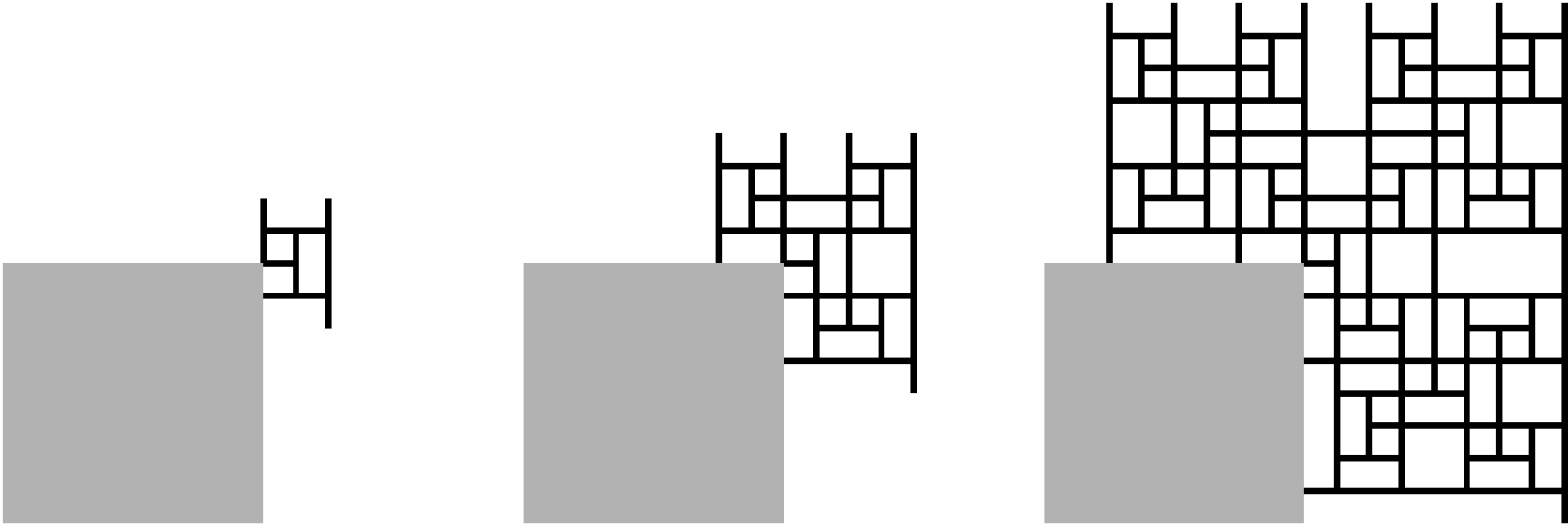}}
\caption{ The corner toothpick structure after $2^k-1$ stages, for $k=2, 3, 4$.}
\label{FigCSUnion}
\end{figure}

An examination of Fig. \ref{FigCS1} and pictures of later stages in
the evolution reveals that after $2^k-1$ stages (for $k \ge 2$)
the structure consists of an essentially solid rectangle of toothpicks 
with one quadrant removed.
The first few cases are shown in Figs. \ref{FigCSUnion}.
More precisely, we have:

\begin{theorem}\label{th1}
After $2^k-1$ stages, for $k \ge 2$, the corner toothpick structure 
is bounded by a rectangle of dimensions 
$(\mbox{height} \times \mbox{width}) =
(2^{k-1}-\frac{1}{2}) \times (2^{k-1}-1)$
with the lower left
$(2^{k-2}-\frac{1}{2}) \times (2^{k-2}-1)$
corner removed and with an additional half-toothpick
protruding downwards from the lower right corner,
in which all the boundary edges are solid rows of
toothpicks except for the top edge which contains no 
horizontal toothpicks,
with a row of $2^{k-1}$ exposed vertical toothpick ends
along the top edge, and
with no exposed toothpick ends in the interior.
Furthermore,
for $k \ge 2$, the number of toothpicks added at the successive stages
while going from stage $2^k$
to stage $2^{k+1}-1$ is given by:
\beql{EqCS1}
c(2^k+i)=
\begin{cases}
2^{k-1}+1,  &\text{if $i=0$}; \\
2c(i) + c(i+1),  &\text{if $i=1, \ldots,2^k-2$}; \\
2c(i) + c(i+1)-1,  &\text{if $i=2^k-1$}. \\
\end{cases}
\eeq
\end{theorem}

\noindent{\bf Proof.}
We use induction on $k$.
The case $k=2$ is readily checked
(cf. Figs. \ref{FigCS1}, \ref{FigCSUnion}).
Suppose the theorem is true for $k$, so that
after the first $2^k-1$ stages we have the structure described
in the theorem. We consider the next $2^k$
stages in the evolution of the bottom right quadrant and the top
two quadrants separately. 

\begin{figure}[htbp]
\begin{center}
\includegraphics[height=2.6in]{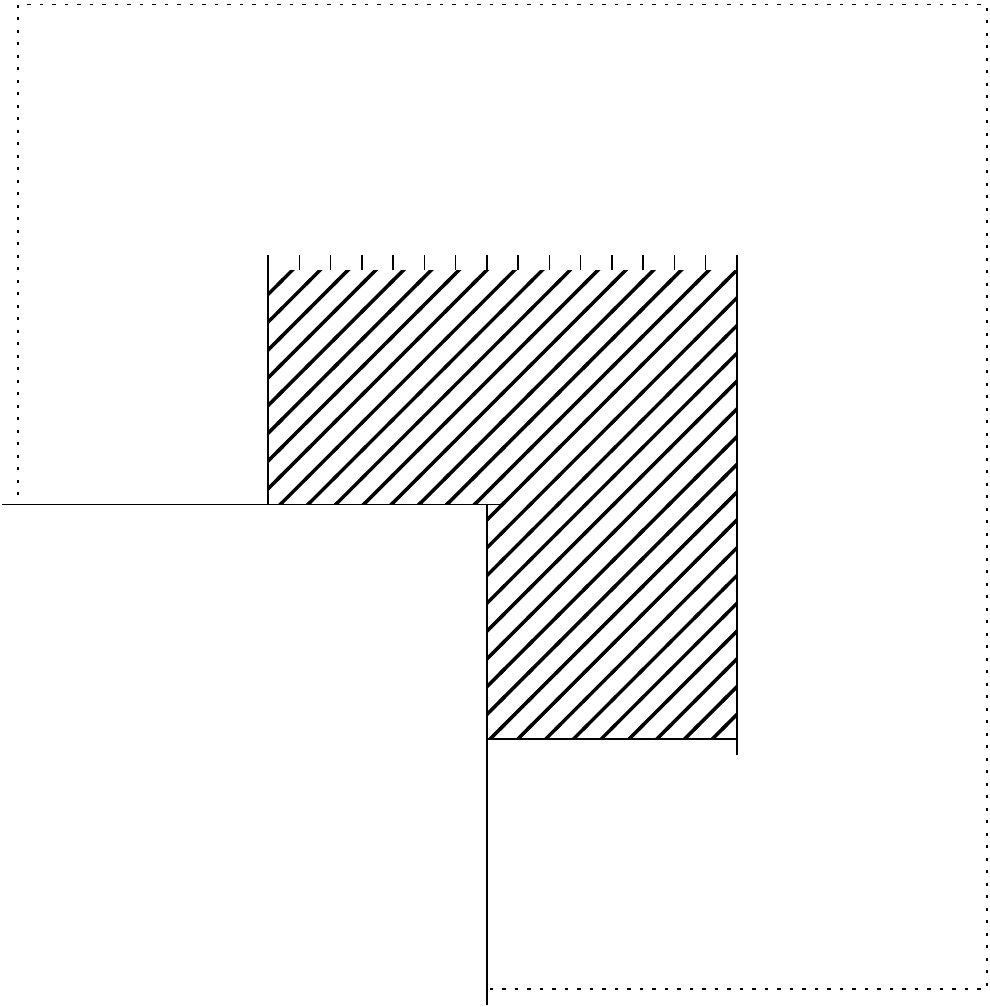}
\end{center}
\caption{Schematic of corner structure after $2^k-1$ stages.}
\label{FigDaveb}
\end{figure}

\begin{figure}[htbp]
\begin{center}
\includegraphics[height=2.6in]{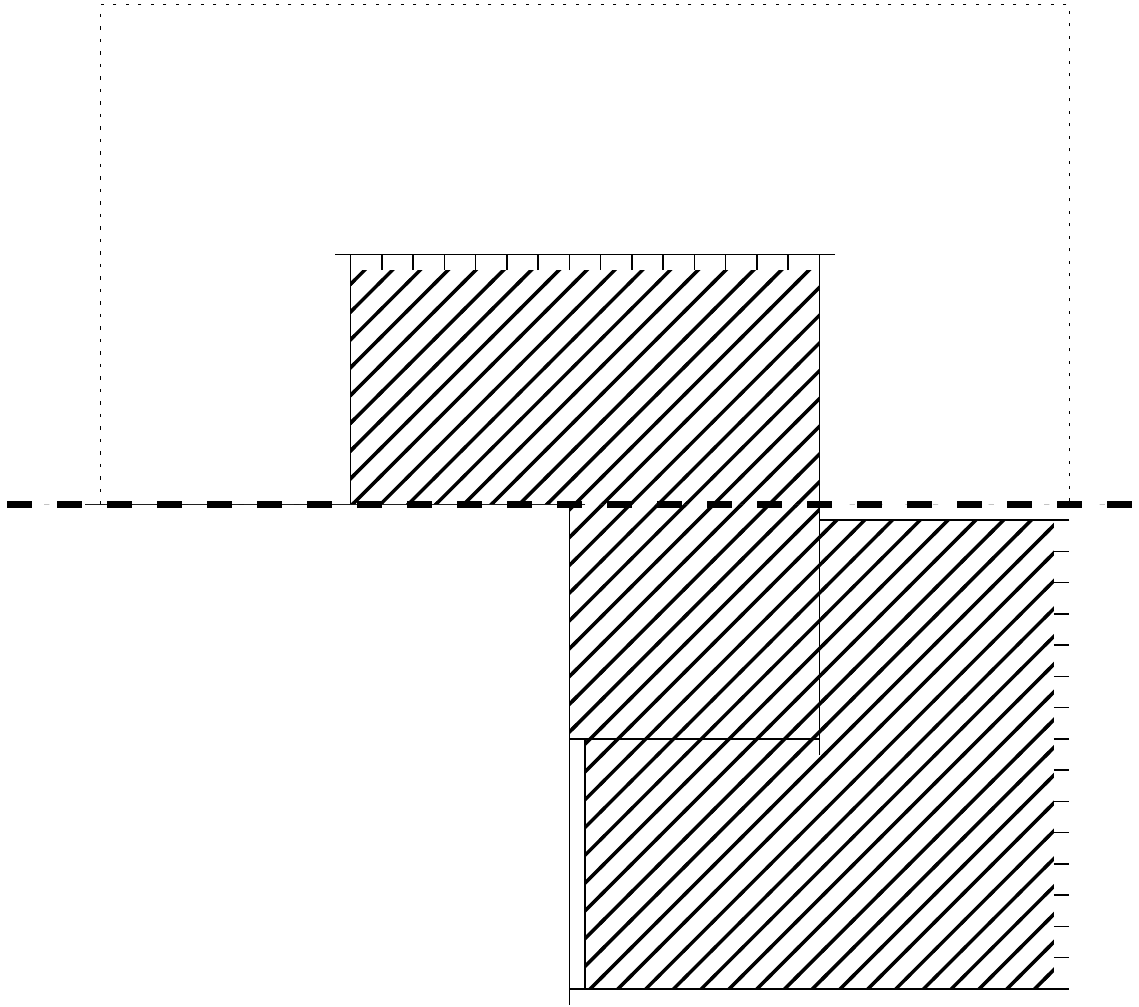}
\end{center}
\caption{Schematic of corner structure after $2^{k+1}-2$ stages
in the third quadrant and $2^k$ stages in the upper half.}
\label{FigDavec}
\end{figure}

\begin{figure}[htbp]
\begin{center}
\includegraphics[height=2.6in]{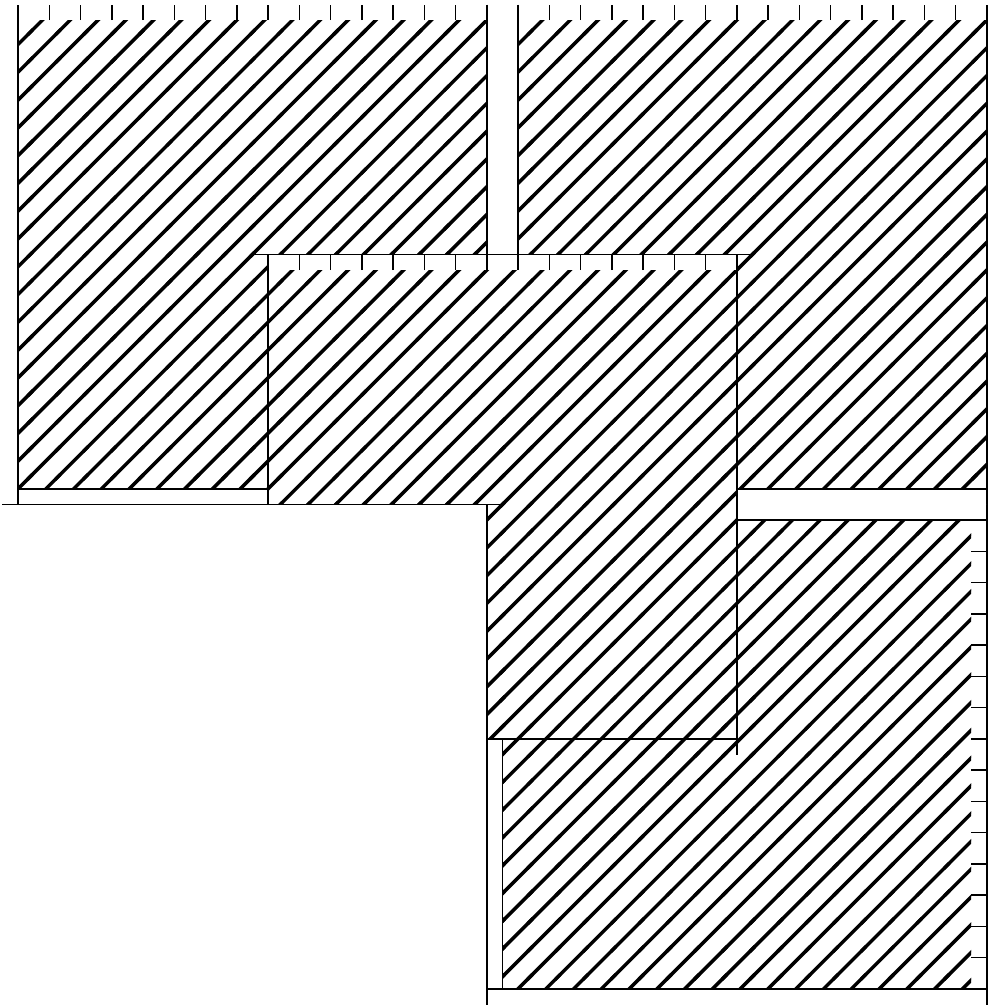}
\end{center}
\caption{Schematic of corner structure after $2^{k+1}-1$ stages.}
\label{FigDaved}
\end{figure}

First, the bottom right quadrant looks like the starting configuration
for the corner structure, with its protruding half-toothpick,
except rotated clockwise by $90^{\circ}$.
So by the induction hypothesis, after further $2^k-1$ stages
we reach a $90^{\circ}$-rotated copy of the $(2^k-1)$-stage structure.
One further step then fills in the right-hand edge,
leaving a half-toothpick protruding downwards from the bottom right corner.
Second, consider what happens to the top half
of the structure. At the first step, the $2^k-1$ vertical
exposed toothpick ends will be covered,
producing overhanging half-toothpicks at the left-
and right-hand ends of the top edge. 
Again these look like the starting configuration
for the corner structure, with the first quadrant a 
mirror image of the second quadrant. So again, by
induction, after a further $2^k-1$ steps we reach the
top half of the desired structure for $k+1$.  
This completes the proof of the first
assertion of the theorem.
(The process is depicted schematically in 
Figs. \ref{FigDaveb}, \ref{FigDavec} and \ref{FigDaved}.)

The recurrence formula \eqn{EqCS1} now follows by 
keeping track of the number of toothpicks that
are added at successive steps as we progress from 
stage $2^k$ to stage $2^{k+1}-1$.~~~$\bsq$

% spot 1

It is worth remarking that this growth in three quadrants, one of which is
a step ahead of the other two, is responsible 
for the terms of the form 
\beql{EqTP0}
2f(i)+f(i+1)
\eeq
which appear in the recurrences
in Theorems  \ref{th1}, \ref{th2} and \ref{thRT1}.

% spot 2

\begin{figure}[htbp]
\begin{center}
\includegraphics[width=4.5in]{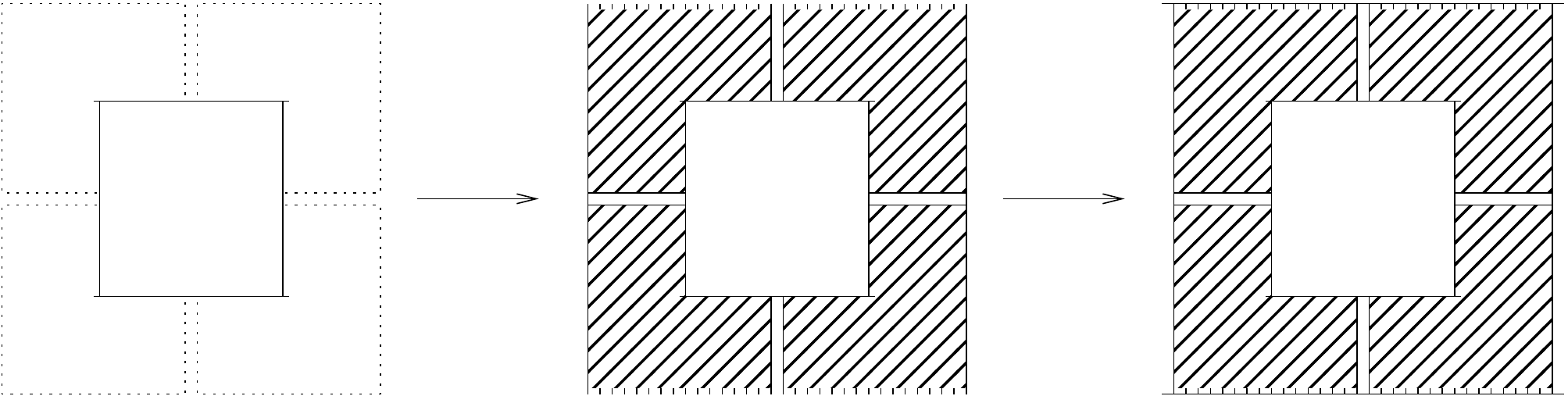}
\end{center}
\caption{Schematic of evolution of toothpick structure:
after $2^k$ stages (left-hand figure),
after a further $2^k-1$ stages (center)
and after one more stage (right-hand figure).}
\label{FigDavee}
\end{figure}

\section{ The toothpick sequences}\label{SecTP2}

Similar recurrences hold for the toothpick sequences $t(n)$ and $T(n)$.

\begin{theorem}\label{th2}
For the toothpick structure discussed in \S\ref{SecTP},
the number of toothpicks added at the $n$th stage
is given by $t(0)=0$, $t(1)=1$, and, for $k \ge 1$,
\beql{EqTP1}
t(2^k+i)=
\begin{cases}
2^k,  &\text{if $i=0$}; \\
2t(i) + t(i+1),  &\text{if $i=1, \ldots,2^k-1$}. \\
\end{cases}
\eeq
\end{theorem}

\noindent{\bf Proof.}
An inductive argument similar to that used in the proof
of Theorem \ref{th1} shows that after $2^k$ steps, for
$k \ge 2$, the toothpick structure 
is bounded by a $2^{k-1} \times (2^{k-1}-1)$
rectangle, with half-toothpicks
protruding horizontally from the four corners,
in which all the boundary edges are solid rows of
toothpicks, and
with no exposed toothpick ends in the interior.
(The cases $k=2$ and $3$ can be seen in Fig. \ref{FigTP1}.)
In the induction step, each quadrant grows 
like a suitably rotated version of the corner structure.
(The evolution is depicted schematically in Fig. \ref{FigDavee}.)
The recurrence formula \eqn{EqTP1} now follows by
keeping track of the number of toothpicks that
are added as we progress from
stage $2^k$ to stage $2^{k+1}-1$.~~~$\bsq$

%\vspace*{+.2in}

\begin{coro}\label{corT}
For the toothpick sequence $T(n)$, we have
$T(0)=0$, and, for $k \ge 0$,
\beql{EqTPT1}
T(2^k+i)=
\begin{cases}
\frac{1}{3}(2^{2k+1}+1)  &\text{if $i=0$}; \\
T(2^k) + 2T(i) + T(i+1) -1,  &\text{if $i=1, \ldots,2^k-1$}. \\
\end{cases}
\eeq
\end{coro}

\noindent{\bf Proof.}
This follows easily from 
$T(n) = \sum_{i=0}^{n} t(i)$,
$$
\sum_{i=2^k}^{2^{k+1}-1} t(i) = 2^k(2^{k+1}-1),
$$
and \eqn{EqTP1}. 
We omit the details.~~~$\bsq$

%%%%%  Triangles

\begin{table}[htbp]
$$
\begin{array}{c|rrrrrrrrrrrrrrrr}
k & \multicolumn{16}{l}{\mbox{~terms~} 2^k, 2^k+1, \ldots, 2^{k+1}-1} \\
%\hline
\cline{1-16}
0   & 0 &   &   &   &   &   &   &   &   &   &   &   &   &   &   &  \\
1   & 1 &   &   &   &   &   &   &   &   &   &   &   &   &   &   &  \\
2   & 2 & 4 &   &   &   &   &   &   &   &   &   &   &   &   &   &  \\
4   & 4 & 4 & 8 & 12 &   &   &   &   &   &   &   &   &   &   &   &  \\
8   & 8 & 4 & 8 & 12 & 12 & 16 & 28 & 32 &   &   &   &   &   &   &   &  \\
16 & 16 & 4 & 8 & 12 & 12 & 16 & 28 & 32 & 20 & 16 & 28 & 36 & 40 & 60 & 88 & 80  \\
\ldots & \multicolumn{16}{l}{\ldots}
\end{array}
$$
\caption{Initial terms of toothpick sequence $t(n)$ arranged in triangular form.}
\label{TabTPtri}
\end{table}

A convenient way to visualize the recurrences \eqn{EqCS1}, \eqn{EqTP1}
and \eqn{EqTPT1}
is to write the sequences $c(n)$, 
$t(n)$ and $T(n)$ as triangular arrays, with $1$, $1$, $2$,
$4$, $8$, $16$, $32$, $\ldots$ terms in the successive rows.
For example, the initial terms of the $t(n)$ sequence are shown in the array
in Table \ref{TabTPtri}.
The row labeled $8$, for instance, begins with $t(8)=8$,
and then, using \eqn{EqTP1} and referring back to the top of the triangle,
continues with the values
$t(9) = 2t(1)+t(2)=4$, 
$t(10) = 2t(2)+t(3)=8$, 
$t(11) = 2t(3)+t(4)=12$, and so on (a kind of ``bootstrap'' process).

%\vspace*{+.2in}

To see a direct connection between the toothpick sequences and the corner
sequences, it is convenient to define $Q(n) := (T(n)-3)/4$ ($n \ge 3$), 
with $Q(0)=Q(1)=0$. This is the number of toothpicks whose centers
are in the interior of the first (or second, third or fourth)
quadrants of the toothpick structure.
Also let $q(n) := Q(n)-Q(n-1)$ ($n\ge 1$) with $q(0)=0$.
The argument used in the proof of Theorem \ref{th1}
shows that
\beql{EqCQ1}
C(n) = 2Q(n)+Q(n+1) + 2, \quad \mbox{~for~} n \ge 2 \,.
\eeq
Hence by taking differences we have $q(n) = t(n)/4$,
\beql{EqCQ2}
c(n) = 2q(n)+q(n+1), \quad \mbox{~for~} n \ge 3 \,,
\eeq
and
\beql{EqCQ2bis}
c(n) = \frac{1}{2} t(n) + \frac{1}{4} t(n+1), \quad \mbox{~for~} n \ge 1 \,.
\eeq

\section{ Rectangles in the toothpick structure}\label{SecWin}
Examination of Figs. \ref{FigTP1}--\ref{Fig64} suggests
that, after any finite number of stages, the toothpick structure
divides the plane into an unbounded region and a number
of squares and rectangles (and no other closed polygonal regions appear).
Let $R(n)$ denote the number of squares and rectangles in
the toothpick structure after $n$ stages, and let $r(n) := R(n)-R(n-1)$
be the number of squares and rectangles that are added at the $n$th stage.
Similarly, let $\rho(n)$ be the number of squares and rectangles that
are added to the corner structure at the $n$th stage.
The initial values of $\rho(n)$, $r(n)$ and $R(n)$ are shown in 
Table \ref{TabRT1} 
(these are entries A168131, A160125 and A160124 in \cite{OEIS}).

\begin{table}[htbp]
$$
\begin{array}{|c|rrrrrrrrrrrrrrrr|}
\hline
n       & 0 & 1 & 2 & 3 & 4 & 5 & 6 & 7 & 8 & 9 & 10 & 11 & 12 & 13 & 14 & 15  \\
\rho(n) & 0 & 0 & 1 & 2 & 1 & 1 & 5 & 7 & 3 & 1 & 4 & 5 & 3 & 7 & 18 & 19 \\
r(n)    & 0 & 0 & 0 & 2 & 2 & 0 & 4 & 10 & 6 & 0 & 4 & 8 & 4 & 4 & 20 & 30 \\
R(n)    & 0 & 0 & 0 & 2 & 4 & 4 & 8 & 18 & 24 & 24 & 28 & 36 & 40 & 44 & 64 & 94 \\
\hline
\end{array}
$$
\caption{The sequences $\rho(n)$, $r(n)$ and $R(n)$ 
for $0 \le n \le 15$.}
\label{TabRT1}
\end{table}

Then an inductive argument, similar to that used to establish
Theorem \ref{th1}, shows the following.

\begin{theorem}\label{thRT1}
All internal regions in the corner and toothpick structures
are squares and rectangles. 
Furthermore, $\rho(0)=\rho(1)=0$, 
$\rho(2)=1$, $\rho(3)=2$, and, for $k \ge 2$,
\beql{EqRT1}
\rho(2^k+i)=
\begin{cases}
2^{k-1}-1,  &\text{if $i=0$}; \\
2\rho (i) + \rho (i+1),  &\text{if $i=1, \ldots,2^k-3$}; \\
2\rho (i) + \rho (i+1)+1,  &\text{if $i=2^k-2$}. \\
2\rho (i) + \rho (i+1)+2,  &\text{if $i=2^k-1$}; \\
\end{cases}
\eeq
and $r(0)=r(1)=r(2)=0$, $r(3)=2$, and, for $k \ge 2$,
\beql{EqRT2}
r(2^k+i)=
\begin{cases}
2^{k}-2,   &\text{if $i=0$}; \\
4\rho (i), &\text{if $i=1, \ldots,2^k-2$}; \\
4\rho (i)+2, &\text{if $i=2^k-1$}. \\
\end{cases}
\eeq
\end{theorem}
We omit the proof.

\section{ The fractal-like structure}\label{SecFrac}
\begin{figure}[htbp]
\begin{center}
\includegraphics[width=4.5in]{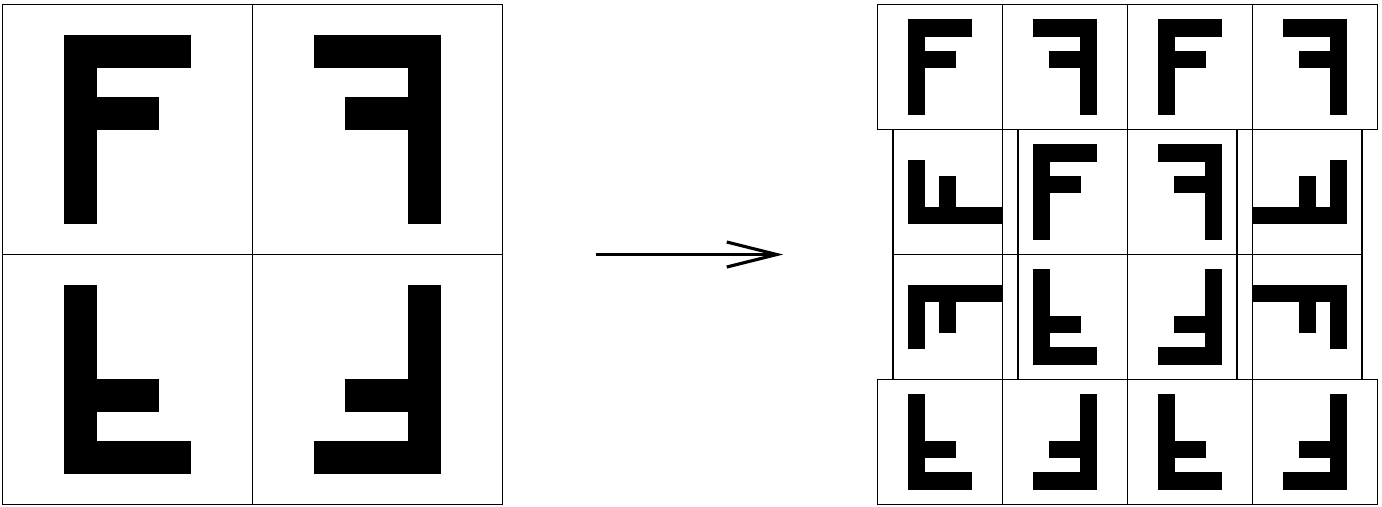}
\end{center}
\caption{Fractal-like transformation recursion, general step.}
\label{FigDavef}
\end{figure}

\begin{figure}[htbp]
\begin{center}
\includegraphics[width=4.5in]{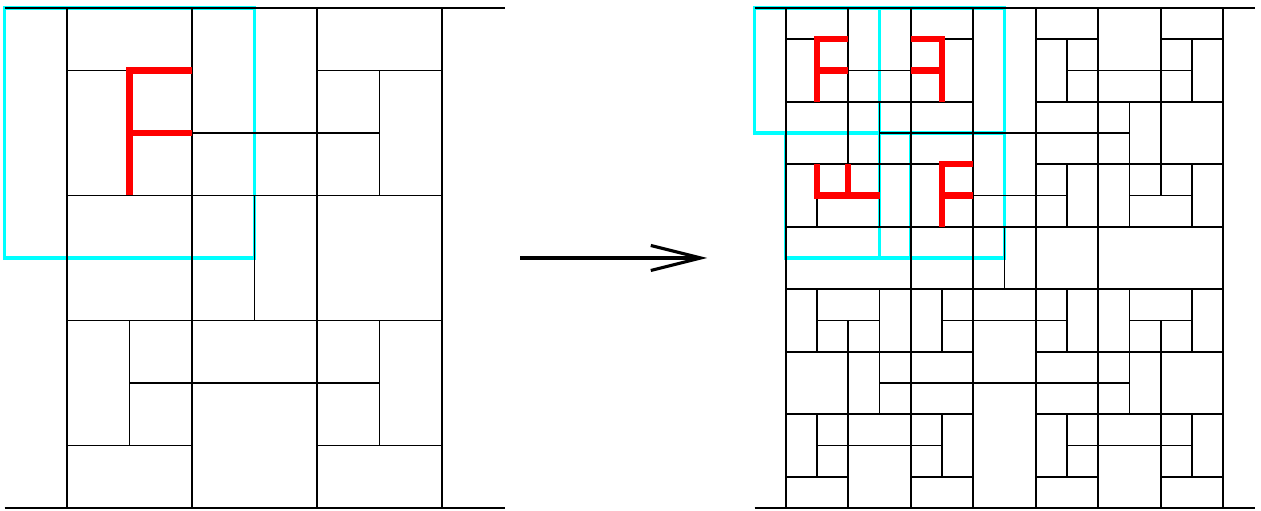}
\end{center}
\caption{Fractal-like recursion going from stage $8$ to stage $16$.}
\label{FigDaveg}
\end{figure}

The recursive structure established in the proofs
of Theorems \ref{th1} and \ref{th2}
also explains the fractal-like appearance
of the toothpick array.
After applying one round of the corner recursion to each
quadrant and then rescaling,
we have the transformation shown schematically in 
Fig. \ref{FigDavef}
(an ``$\F$'' is used to indicate orientation of the
various pieces), and in a specific example in Fig.~\ref{FigDaveg}.  
Note that four of the blocks (the sideways ``$\F$''s) are
shifted by a half-toothpick towards the center.
Because of this shift the toothpick structure
is not strictly self-similar (cf. \cite{Falc90})
and so is not a true fractal. 
The same is true for all the structures
we will meet in  this paper: they have a fractal-like
growth, but are not strictly self-similar.

%%%%% FIG 12 %%%%%%%%%%%%
\begin{figure}[htbp]
\begin{center}
%\includegraphics[ width=4.3in, trim=0in 0.2in 0.4in 0.8in ]{jubin6.pdf}
% Using
% \includegraphics[ width=3.3in, trim=0.15in 0.3in 0.55in 0.85in ]{jubin6.pdf}
% and
% \includegraphics[ width=3.3in, trim=0.15in 0.3in 0.55in 0.85in ]{jubin7.pdf}
% they both fit on the same page (or widths below 3.3in).
% For width=3.35in, the figures show up on separate pages, but with other
% text.  I stumbled on this by accident, but think it looks pretty good.
% For width=3.4in or larger, they end up by themselves on separate pages.
\includegraphics[ width=3.3in, trim=0.15in 0.3in 0.55in 0.85in ]{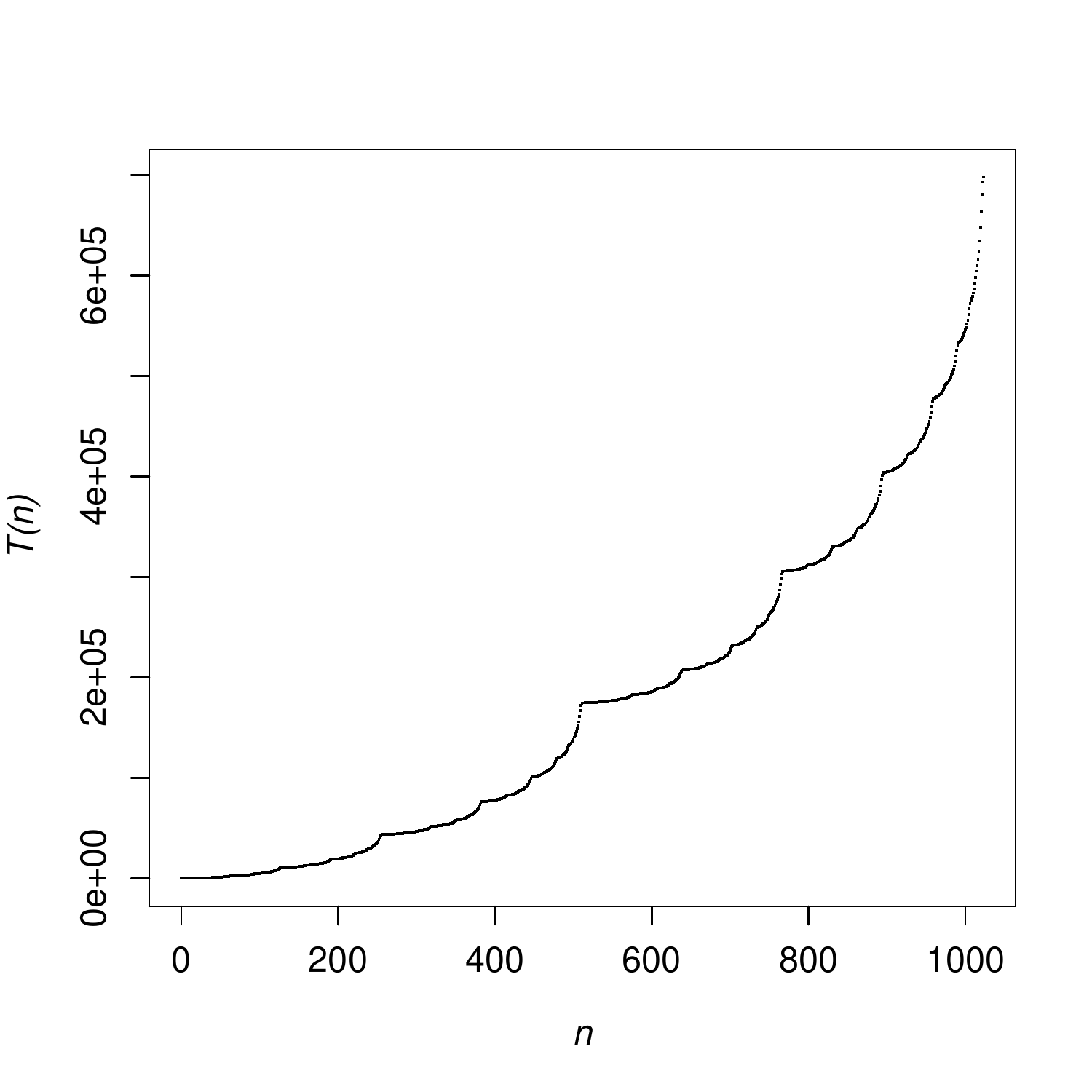}
\end{center}
\caption{Plot of toothpick sequence $T(0), \ldots, T(1023)$.}
\label{FigJubin}
\end{figure}

\begin{figure}[htbp]
\begin{center}
\includegraphics[ width=3.3in, trim=0.15in 0.3in 0.55in 0.85in ]{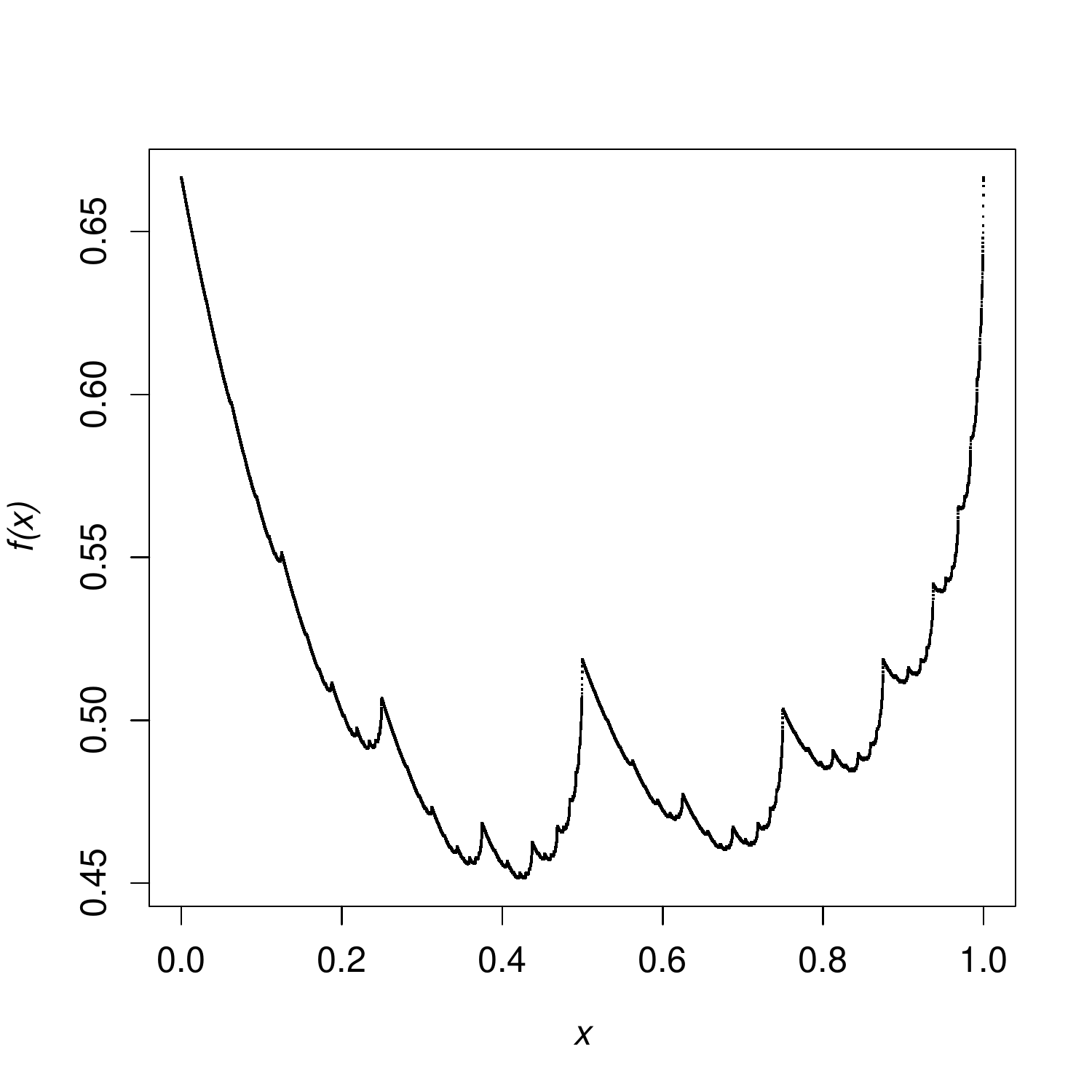}
\end{center}
\caption{
Plot of 
$(x=\frac{i}{2^k}, f(x)=\frac{T(n)}{n^2})$, $n=2^k+i$, $0 \le i <2^k$,
for $k=14$.}
\label{FigJubin2}
\end{figure}

A plot of $T(0), \ldots, T(2^k-1)$ for increasing
values of $k$ shows that the points
lie roughly on a parabola, with irregularities
caused by the fractal-like behavior (see for example Fig. \ref{FigJubin}).
%can be seen by using the ``graph'' button
%in entry A139250 in \cite{OEIS}).
Beno\^{i}t Jubin \cite{Jubi10} has investigated
$\varlimsup T(n)/n^2$ and $\varliminf T(n)/n^2$. His
numerical results suggest that 
$\varlimsup T(n)/n^2 = 2/3$, with local maxima at 
the values $n = 2^k-1$, and 
$\varliminf T(n)/n^2 \approx 0.4513058$ with local minima at 
the following values of $n$:
$$
1, 2, 5, 12, 21, 44, 89, 180, 362, 728, 1459, 2921, \ldots 
$$
(see A170927 for further terms).
His upper limit can be established from Corollary \ref{corT}:

\begin{theorem}\label{thBJ}
For $n \ge 1$, 
\beql{EqBJ1}
\frac{T(n)}{n^2} \le \frac{2}{3} + \frac{1}{3n}\,,
\eeq
with equality if and only if $n=2^k-1$ for some $k$.
Hence $\varlimsup T(n)/n^2 = 2/3$.
\end{theorem}

\noindent{\bf Proof.}
The result is true for $n \le 3$,
and for $n >3$ for the special values
$n=2^k-1$, when
$T(n)=(2^k-1)(2^{k+1}-1)/3$, so
$$
\frac{T(n)}{n^2} = \frac{2}{3} + \frac{1}{3n}\,,
$$
and 
$n=2^k$, when
$T(n)=(2^{2k+1}+1)/3$, so
$$
\frac{T(n)}{n^2} = \frac{2}{3} + \frac{1}{3n^2}
< \frac{2}{3} + \frac{1}{3n}\,.
$$
For the general case we use induction,
and assume that \eqn{EqBJ1} holds for all $n \le 2^k$.
Let $n=2^k+i$, $1 \le i \le 2^k-2$.
Then \eqn{EqBJ1} follows from \eqn{EqTPT1}
and the induction hypothesis.~~~$\bsq$

Jubin also observes that there is a continuous
%but nowhere differentiable 
function on $[0,1]$ that describes
the asymptotic behavior of $T(n)$.
This is the function whose graph is the Hausdorff limit
of the finite sets $E_k$ consisting of the points
$(x=\frac{i}{2^k}, f(x)=\frac{T(n)}{n^2})$ for $n=2^k+i$, $0 \le i <2^k$.
This function takes the value $\frac{2}{3}$ at 
$x=0$ and $x=1$, and has its minimum at around
(0.427451, 0.4513058). It is non-differentiable
at the dyadic rational points between $0$ and $1$. 
Figure \ref{FigJubin2} shows $E_{14}$.

% FIG 14
\begin{figure}[htbp]
%\centerline{\includegraphics[width=5in]{UW.pdf}}
\centerline{\includegraphics[width=4.1in]{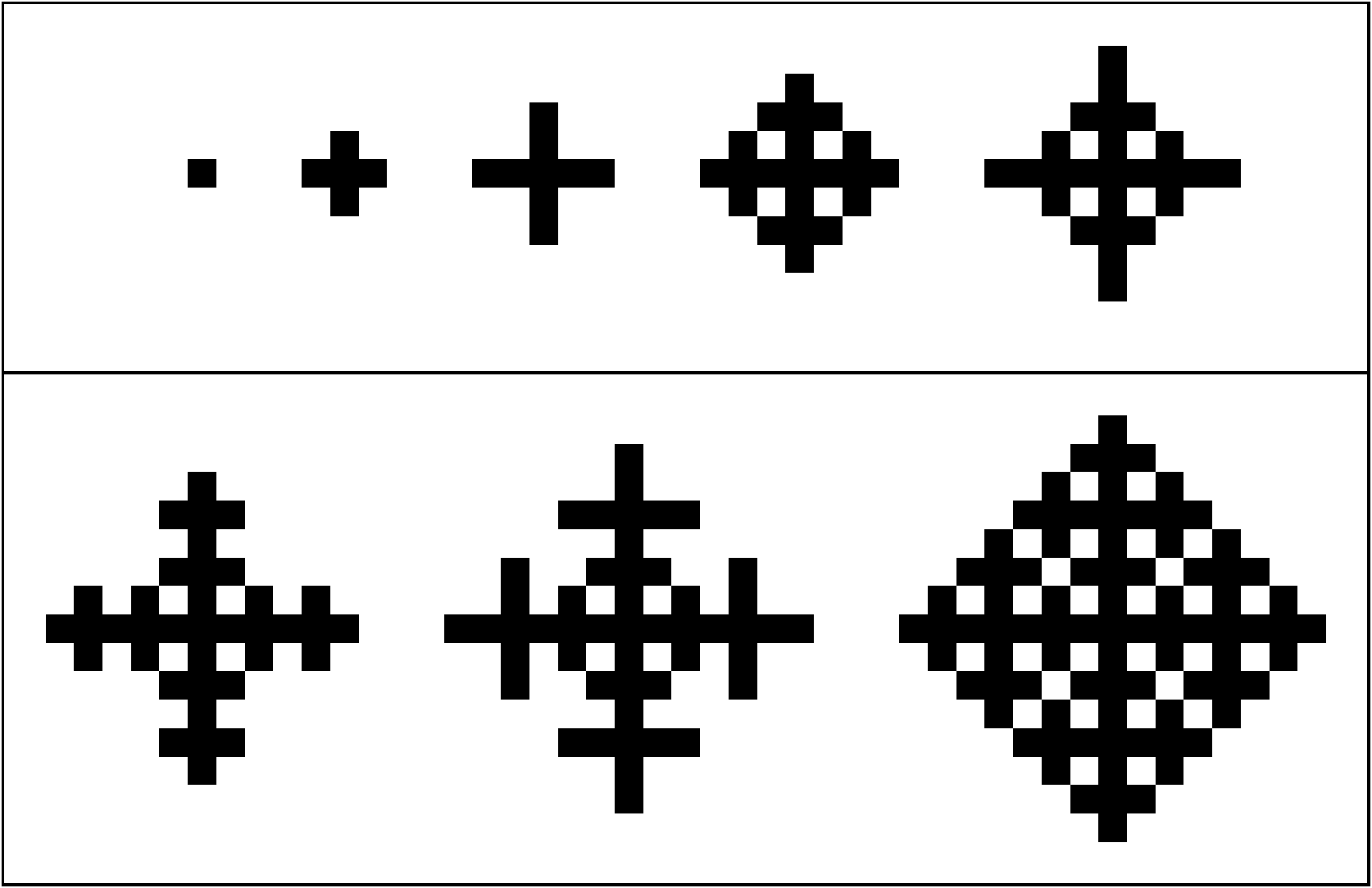}}
\caption{Stages $1$ through $8$ of the
evolution of the Ulam-Warburton structure.
The numbers of $\ON$ cells in the successive stages, $U(1), \ldots,
U(8)$, are $1, 5, 9, 21, 25, 37, 49, 85$.}
\label{FigUW1}
\end{figure}

%%%%%%%%%%%%%%%%%%%%%%%%%%%%%%%%%%%%%%
% Section: 
%%%%%%%%%%%%%%%%%%%%%%%%%%%%%%%%%%%%%%

\section{The Ulam-Warburton cellular automaton}\label{SecUW}

As we will see in \S\ref{SecCAG},
the toothpick structure can be modeled by a cellular automaton
on a planar graph.  In this section, we consider
a simpler example of the same type, the
arrangement of square cells generated by the Ulam-Warburton
cellular automaton (Ulam \cite{Ulam62}, Singmaster \cite{Sing03},
Stanley and Chapman \cite{Stan94}, Wolfram \cite[p.~928]{NKS}).
The cells are the squares in an infinite square grid, and the neighbors 
of each cell are defined to be the four squares which
share an edge with it.
(This is the von Neumann neighborhood of the cell, in
the notation of \cite{Mitch96}.)
At stage $0$, no cells are $\ON$.
At stage $1$, a single cell is turned $\ON$.
Thereafter, a cell is changed from $\OFF$ to $\ON$ at stage $n$ if and only if
exactly one of its four neighbors was $\ON$ at stage $n-1$.
Once a cell is $\ON$ it stays $\ON$.
This is ``Rule 686'' in the notation of \cite{PaWo85}, \cite{NKS}.

Let $u(n)$ ($n\ge 0$) denote the number of
cells that are changed from $\OFF$ to $\ON$ at the $n$th stage,
and let $U(n) := \sum_{i = 0}^{n} u(i)$ be the total number
of $\ON$ cells after $n$ stages.
The initial values of $u(n)$ and $U(n)$ are shown in Table \ref{TabUW1}.
These sequences are respectively entries A147582 and A147562 in \cite{OEIS}.
Fig. \ref{FigUW1} shows stages $1$ through $8$ of the
evolution of the this structure.
(As is suggested by Fig. \ref{FigUW1} and more
particularly by the movie 
linked to entry A147562, this structure also has a fractal-like growth.)

\begin{table}[htb]
$$
\begin{array}{|c|rrrrrrrrrr|}
\hline
n   & 0 & 1 & 2 & 3 & 4 & 5 & 6 & 7 & 8 & 9 \\
u(n) & 0 & 1 & 4 & 4 & 12 & 4 & 12 & 12 & 36 & 4 \\
U(n) & 0 & 1 & 5 & 9 & 21 & 25 & 37 & 49 & 85 & 89 \\
\hline
n & 10 & 11 & 12 & 13 & 14 & 15 & 16 & 17 & 18 & 19 \\
u(n) & 12 & 12 & 36 & 12 & 36 & 36 & 108 & 4 & 12 & 12 \\
U(n) & 101 & 113 & 149 & 161 & 197 & 233 & 341 & 345 & 357 & 369 \\
\hline
n & 20 & 21 & 22 & 23 & 24 & 25 & 26 & 27 & 28 & 29 \\
u(n) & 36 & 12 & 36 & 36 & 108 & 12 &  36 & 36 & 108 & 36 \\
U(n) & 405 & 417 &  453 & 489 & 597 & 609 & 645 & 681 & 789 & 825 \\
\hline
n & 30 & 31 & 32 & 33 & 34 & 35 & 36 & 37 & 38 & 39 \\
u(n) & 108 & 108 & 324 & 4 & 12 & 12 & 36 & 12 & 36 & 36  \\
U(n) & 933 & 1041 & 1365 & 1369 & 1381 & 1393 & 1429 & 1441 &  1477 & 1513 \\
\hline
n & 40 & 41 & 42 & 43 & 44 & 45 & 46 & 47 & 48 & 49 \\
u(n) & 108 & 12 & 36 & 36 & 108 & 36 & 108 & 108 & 324 & 12 \\
U(n) & 1621 & 1633 & 1669 & 1705 & 1813 & 1849 & 1957 & 2065 & 2389& 2401 \\
\hline
\end{array}
$$
\caption{The sequences $u(n)$ and $U(n)$ from
the Ulam-Warburton cellular automaton, for $0 \le n \le 49$.}
\label{TabUW1}
\end{table}

\begin{theorem}\label{th3}
(i) The number of cells that turn from $\OFF$
to $\ON$ at stage $n$ of the Ulam-Warburton cellular automaton
satisfies the recurrence $u(0)=1$, $u(1)=1$,
and, for $k \ge 0$,
\beql{EqUW1}
u(2^k+1+i)=
\begin{cases}
4,  &\text{if $i=0$}; \\
3u(i),  &\text{if $i=1, \ldots,2^k-1$}. \\
\end{cases}
\eeq
(ii) There is an explicit formula: $u(0)=0$, $u(1)=1$
and
\beql{EqUW2}
u(n) = 4 \mydot 3^{\wt(n-1)}~-~1, ~~n\ge 2\,,
\eeq
where $\wt(n)$, the ``binary weight'' of $n$, is the number of $1$'s in
the binary expansion of $n$ (entry A000120 in \cite{OEIS}). \\
(iii) The $u(n)$ have generating function
\beql{EqUW3}
%x + \frac{4x}{3} \left\{ \prod_{k \ge 0} (1+3x^{2^k}) - 1 \right\} \,.
x( -\frac{1}{3} + \frac{4}{3} \prod_{k \ge 0} (1+3x^{2^k}) ) \,.
\eeq
\end{theorem}

\noindent{\bf Proof.}
Part (i) follows by an inductive argument similar to
that used in the proofs of Theorems \ref{th1} and \ref{th2}.
The appropriate ``corner sequence'' is A048883, in which
the $n$th term is $3^{\wt(n-1)}$ ($n\ge 1$),
with partial sums given by A130665.
Part (ii) follows from (i) by induction on $n$.
Part (iii) follows from the generating function for A048883, which is
$\prod_{k \ge 0} (1+3 x^{2^k})$.~~~$\bsq$

\paragraph{Remarks.}

1. Now the corner sequence has three quadrants that grow
in synchronism, so the $2c(i)+c(i+1)$ terms in \eqn{EqCS1} 
are replaced by the $3u(i)$ term in \eqn{EqUW1}.
 
2. Parts (i) and (ii) of the theorem can be found in
Singmaster \cite{Sing03}
and Stanley and Chapman \cite{Stan94}, and
part (i) at least was probably known to
J. C. Holladay and Ulam.  On page 216 of \cite{Ulam62},
Ulam remarks that for certain structures similar to this one
(exactly which ones is left unspecified), J. C. Holladay
showed that ``at generations whose index number $n$ is
of the form $n=2^k$, the growth is cut off everywhere 
except on the `stems', i.e. the straight lines issuing from
the original point.'' This is certainly consistent 
with the recurrence \eqn{EqUW1}.   

3. Two properties of the Ulam-Warburton structure given
in \cite{Stan94} are worth mentioning here. (i) When
considered as a subgraph of the infinite square grid, 
the structure is a tree. This is also
true for the toothpick structure---see \S\ref{SecCAG}.
(ii) Let $n-1 = \sum_{i=1}^{w} 2^{r_i}$
($r_1 > r_2 > \cdots > r_w \ge 0$) be the binary expansion of $n-1$.
Then a necessary and sufficient condition for the cell
at $P=(x,y) \in \ZZ \times \ZZ$ to be turned from $\OFF$ to $\ON$
at stage $n>1$ is that $P = \sum_{i=1}^{w} 2^{r_i} v_i$,
where $v_i \in 
\{(-1,0), (1,0), (0,-1), (0,1)\}$,
subject to $v_i \ne -v_{i-1}$ for $i>1$.
We have no such characterization of the toothpicks added 
at the $n$th stage.
It is a consequence of this (although it is not
mentioned in \cite{Stan94}) that the cells that are
turned $\ON$ at some stage are the cells $(x,y)$ with $x=0$ or
$y=0$, and the cells with $xy \ne 0$ for which the highest
power of $2$ dividing $x$ is different from the highest power of 
$2$ dividing $y$. Again we know of no analog for the toothpick structure.

%%%%%%%%%%%%%%%%%%%%%%%%%%%%%%%%%%%%%%
% Section: 
%%%%%%%%%%%%%%%%%%%%%%%%%%%%%%%%%%%%%%

\section{Leftist toothpicks}\label{SecLT}
Stimulated by Theorem \ref{th3},
we set out to look for analogues of \eqn{EqUW2} and \ref{EqUW3}
for the toothpick sequence $t(u)$.
Our first attempt was a failure, but led to an interesting
connection with Sierpi\'{n}ski's triangle.

\begin{figure}[htbp]
\centerline{\includegraphics[width=1.5in]{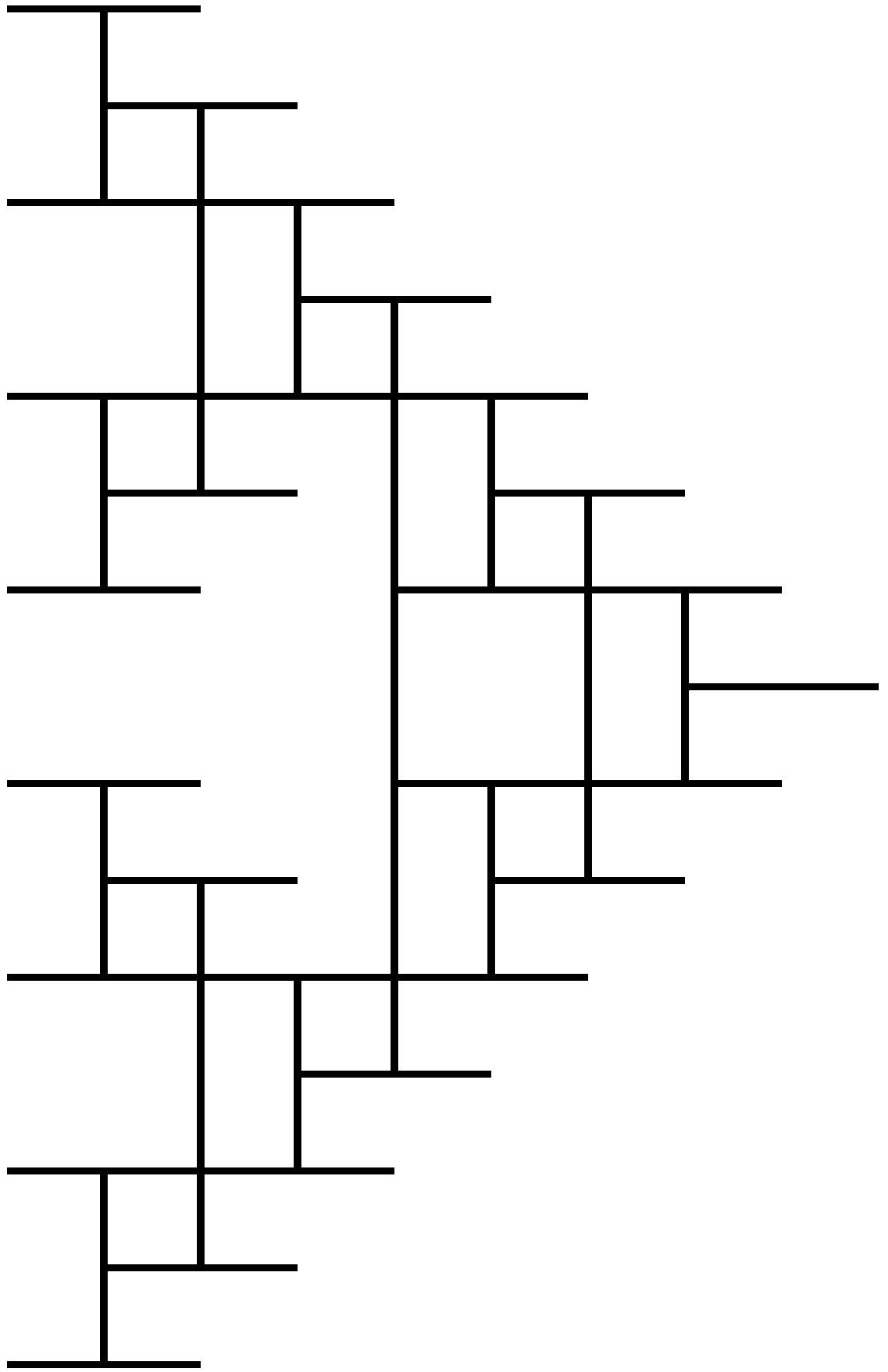}}
\caption{Stages $1$ through $15$ of the
evolution of the leftist toothpick structure.}
\label{FigLT1}
\end{figure}

We define the ``leftist toothpick'' structure as follows.
We start with a single horizontal toothpick at stage $1$,
and extend the structure using the toothpick rule of 
\S\ref{SecTP}, except that if a toothpick is horizontal,
a new toothpick can be added only at its left-hand end.
Let $l(n)$ ($n\ge 1$) denote the number of
toothpicks added at the $n$th stage, with $l(0)=0$,
and let $L(n) := \sum_{i = 0}^{n} l(i)$ be the total number
of toothpicks after $n$ stages.
These are respectively entries A151565 and A151566 in \cite{OEIS}.
The initial values of $l(n)$ and $L(n)$ are shown in Table \ref{TabLT1}.
Figure \ref{FigLT1} shows the first $15$ stages of the evolution
(the starting toothpick is the apex of the triangle, at the right).

\begin{table}[htbp]
$$
\begin{array}{|c|rrrrrrrrrrrrrrrr|}
\hline
n    & 0 & 1 & 2 & 3 & 4 & 5 & 6 & 7 & 8 & 9 & 10 & 11 & 12 & 13 & 14 & 15  \\
l(n) & 0 & 1 & 1 & 2 & 2 & 2 & 2 & 4 & 4 & 2 & 2 &   4 & 4 & 4 & 4 & 8 \\
L(n) & 0 & 1 & 2 & 4 & 6 & 8 & 10 & 14 & 18 & 20 & 22 & 26 & 30 & 34 & 38 & 46 \\
\hline
\end{array}
$$
\caption{The leftist toothpick sequences $l(n)$ and $L(n)$ for $0 \le n \le 15$.}
\label{TabLT1}
\end{table}

The reason for investigating this structure is that,
at least for the early stages, 
%the toothpick structure of \S\ref{SecTP} 
%is made up of four rotated copies of the leftist structure. 
the part of the toothpick structure of \S\ref{SecTP} 
in the $90^{\circ}$ sector $x-2 \le y \le 2-x$
is essentially equal to the leftist structure.
This breaks down, however, at stage $14$. 
Nevertheless, the leftist structure has some interest.  
For if we rotate the structure anticlockwise by $90^{\circ}$ 
and erase all the horizontal toothpicks, we obtain a triangle
in which the vertical toothpicks correspond
exactly to the positions of the $1$s in 
Sierpi\'{n}ski's triangle (i.e., Pascal's triangle read modulo $2$
\cite{Finc03}, \cite{PJS92}, \cite{Schr91}, \cite[Chap. 3]{NKS}),
and the gaps between the vertical toothpicks to the $0$s.
Once observed, this is easy to prove.

Gould's sequence, entry A001316 in \cite{OEIS},
gives the number of odd entries in row $n$ of Pascal's triangle,
which is $2^{\wt(n)}$, with generating function
$\prod_{k \ge 0} (1+2x^{2^k})$.
Allowing for the different offset, we conclude that 
the leftist toothpick sequence $l(n)$ is given by
$l(2n-1)=l(2n)=2^{\wt(n-1)}$.

%%%%%%%%%%%%%%%%%%%%%%%%%%%%%%%%%%%%%%
% Section: 
%%%%%%%%%%%%%%%%%%%%%%%%%%%%%%%%%%%%%%

\section{Generating functions}\label{SecGF}
The following theorem suggests why generating functions of the form
\eqn{EqGF1} arise in connection with recurrences
of the form \eqn{EqCS1}, \eqn{EqTP1}.

\begin{theorem}\label{th4}
Given integers $\al, \be, \ga, \de$, let
the Taylor series expansion of
\beql{EqGF2}
x(\al + \be x) \prod_{k \ge 1} (1 + \ga x^{2^k - 1} + \de x^{2^k}) 
\eeq
be $a(0) + a(1)x + a(2)x^2 + \cdots$. Then we have
$a(0) = 0, a(1)=\al,$
%  \quad a(1)=\al\ga+\be, \quad a(2)=\al\de+\be\ga, \quad a(3)=\al\ga+\be\de,
%a(0)=\al,  \quad a(1)=\al\ga+\be, \quad a(2)=\al\de+\be\ga, 
and for $n\ge 2$, 
\beql{EqGF3b}
a(2^k+i)=
\begin{cases}
\al\ga + \be{\de}^{k-1}  &\text{if $i=0$}; \\
\de a(i) + \ga a(i+1),  &\text{if $i=1, \ldots,2^k-2$}; \\
\de a(i) + \ga a(i+1) -\al{\ga}^2,  &\text{if $i=2^k-1$}. \\
\end{cases}
\eeq
\end{theorem}

\noindent{\bf Proof.}
The generating function is
$$
x(\al + \be x) 
(1 + \ga x + \de x^{2}) 
(1 + \ga x^{3} + \de x^{4}) 
(1 + \ga x^{7} + \de x^{8}) 
\cdots .
$$
Consider the coefficient of $x^{10}$, say.
The only way to build up $x^{10}$ is
to combine the term $\de x^8$ with $a(2)x^2$,
or 
the term $\ga x^7$ with $a(3)x^3$.
Hence $a(2^3+2)=\de a(2) + \ga a(3)$.
Similar arguments holds for the general case,
although adjustments are needed when
$i=0$ or $2^k-1$. We omit the details.~~~$\bsq$

An analogous result holds if the product in \eqn{EqGF2}
starts at $k=0$.

Notice the resemblance between Equations \eqn{EqCS1}, \eqn{EqTP1}
and \eqn{EqGF3b}.  In particular, we have:

\begin{theorem}\label{th5}
The generating function for the corner sequence $c(n)$ is
\beql{EqGCS1}
\sC (x) := x(1+x) \prod_{k \ge 1} (1 + x^{2^k - 1} + 2 x^{2^k}) \, .
\eeq
The generating function for the toothpick sequence $t(n)$ is
\begin{align}\label{EqGTP1}
   \sT (x) & :=  ~ \frac{x}{1+2x} \{
1 +4x(1+x) \prod_{k \ge 1} (1 + x^{2^k - 1} + 2 x^{2^k}) \}  \nonumber \\
        {} & =  ~ \frac{x}{1+2x} \{
1 +2x \prod_{k \ge 0} (1 + x^{2^k - 1} + 2 x^{2^k}) \}  \,,
\end{align}
and therefore the generating function for the toothpick sequence $T(n)$ is
\beql{EqGTP4}
\frac{x}{(1-x)(1+2x)} \{
1 +2x \prod_{k \ge 0} (1 + x^{2^k - 1} + 2 x^{2^k}) \}  \,.
\eeq
\end{theorem}

\noindent{\bf Proof.}
For the first assertion, we set $\al = \be = \ga = 1$, $\de =2$
in Theorem \ref{th4} and use \eqn{EqCS1}.
For the second assertion, we note that \eqn{EqCQ2} implies
that the generating functions for $c(n)$ and $q(n)$
are related by
\beql{EqGTP2}
\sC (x) = x + x^2 + (2 + \frac{1}{x}) \sQ (x) \,,
\eeq
and by definition we have
\beql{EqGTP3}
\sT (x) = x + 2x^2 + \sQ (x) \,.
\eeq
Eliminating $\sQ (x)$, we obtain \eqn{EqGTP1}.~~~$\bsq$

\paragraph{Remark.}
Equation \eqn{EqGTP4} was conjectured by Gary W. Adamson \cite{Adam09}.
Consider the sequence 
$1,1,2,1,3,4,4,1,3,4,5,\ldots$ 
with generating function
$$
\prod_{k \ge 1} (1 + x^{2^k - 1} + 2 x^{2^k})
$$
(entry A151550 in \cite{OEIS}).
Adamson discovered that if this sequence is convolved
with the sequence $1,2,2,2,2,2,\ldots$, the result
appeared to coincide with the corner sequence $C(n)$.
When expressed in terms of generating functions, his
conjecture is essentially equivalent to \eqn{EqGTP4}.

%%%%%%%%%%%%%%%%%%%%%%%%%%%%%%%%%%%%%%
% Section: 
%%%%%%%%%%%%%%%%%%%%%%%%%%%%%%%%%%%%%%

\section{Explicit formulas}\label{SecEX}
A second comment in \cite{OEIS},
this time from Hagen von Eitzen, was instrumental in the
discovery of explicit formulas for many of these sequences.
Von Eitzen \cite{Eitz09} contributed the sequence
$2,3,3,3,5,6,4,3,5,6,6,\ldots$ 
with generating function 
$\prod_{k \ge 0} (1 + x^{2^k - 1} + x^{2^k})$
to \cite{OEIS} (it is entry A160573)
and provided an elegant explicit formula for the $n$th term:
\beql{EqHVE1}
\sum_{m \ge 0} \binom{\wt(n+m)}{m} \,.
\eeq

This can be generalized. For this it is convenient
to omit the initial linear factors from \eqn{EqGF1}
but to start the product at $k=0$.

\begin{theorem}\label{th6}
Let the Taylor series expansion of
\beql{EqHVE2}
\prod_{k \ge 0} (1 + \ga x^{2^k - 1} + \de x^{2^k}) 
\eeq
be $a(0) + a(1)x + a(2)x^2 + \cdots$. Then 
\beql{EqHVE3}
a(n) =
\sum_{m \ge 0} {\ga}^{m} {\de}^{\wt(n+m)-m} \binom{\wt(n+m)}{m} \,.
\eeq 
\end{theorem}

\noindent{\bf Proof.} (Based on von Eitzen's proof of \eqn{EqHVE1}.)
First, observe that 
\beql{EqHVE4}
\prod_{k \ge 0} (1 + \de x^{2^k}) =
\sum_{n = 0}^{\infty} {\de}^{\wt(n)} x^n \,,
\eeq 
since when getting a term $x^n$, we pick up a factor of $\de$
for every $1$ in the binary expansion of $n$.
When we expand
\beql{EqHVE5}
\prod_{k \ge 0} (1 + \ga x^{2^k - 1} + \de x^{2^k}) \,,
\eeq
instead of the product in \eqn{EqHVE4},
each time we replace a term $\de x^{2^k}$ by $\ga x^{2^{k-1}}$,
we lose a factor of $x$ in the product, but we gain because there 
may be several ways to choose the factors in which to do the replacement.
Suppose we do this replacement in $m$ of the terms in
\eqn{EqHVE5}. Then we must increase $n$ to $n+m$,
we gain by a factor of $\binom{\wt(n+m)}{m}$, but we have 
to replace $m$ factors of $\de$ by $\ga$s, for a net contribution of
${\ga}^{m} {\de}^{\wt(n+m)-m} \binom{\wt(n+m)}{m}$ to the sum.~~~$\bsq$

\paragraph{Remark.}
Note that there are only finitely many nonzero terms in the summations
\eqn{EqHVE1} and \eqn{EqHVE3}.  For large $n$ the number
of nonzero terms is roughly $\log_{2} n$. More precisely,
the number of nonzero terms for any $n$
is given by entry A100661 in \cite{OEIS}.

We can use Theorem \ref{th6} to obtain an explicit formula
for the toothpick sequence $t(n)$.

\begin{theorem}\label{th7}
\beql{EqHVE6}
t(2^k+1+i)=
\begin{cases}
2 \sum_{m \ge 0} {2}^{\wt(i+m)-m} \binom{\wt(i+m)}{m},
   &\text{if $0 \le i \le 2^k-2$}; \\
2^{k+1},  &\text{if $i=2^k-1$}. \\
\end{cases}
\eeq
\end{theorem}

\noindent{\bf Proof.} 
This theorem is an instance where our decision to start enumerations
at $n=0$ rather than $n=1$ (as discussed in Remark $1$ 
in \S\ref{SecTP}) causes complications. The result would be
more elegant if we had made the other choice.
So, just for this proof, let us define $\hat{t}(n) = t(n+1)$
for $n \ge 0$.  Then from Theorem \ref{th2},
$\hat{t}(n)$ satisfies the recurrence
\beql{EqTP1bis}
\hat{t}(2^k+i)=
\begin{cases}
2\hat{t}(i) + \hat{t}(i+1),  &\text{if $i=0, \ldots,2^k-2$}; \\
2^{k+1},  &\text{if $i=2^k-1$}. \\
\end{cases}
\eeq

\begin{table}[htbp]
$$
\begin{array}{c|rrrrrrrrrrrrrrrr}
k & \multicolumn{16}{l}{\mbox{~terms~} 2^k, 2^k+1, \ldots, 2^{k+1}-1} \\
%\hline
\cline{1-16}
0   & 1 &   &   &   &   &   &   &   &   &   &   &   &   &   &   &  \\
1   & 2 &   &   &   &   &   &   &   &   &   &   &   &   &   &   &  \\
2   & 4 & 4 &   &   &   &   &   &   &   &   &   &   &   &   &   &  \\
4   & 4 & 8 & 12 & 8 &   &   &   &   &   &   &   &   &   &   &   &  \\
8   & 4 & 8 & 12 & 12 & 16 & 28 & 32 & 16 &   &   &   &   &   &   &   &  \\
16  & 4 & 8 & 12 & 12 & 16 & 28 & 32 & 20 & 16 & 28 & 36 & 40 & 60 & 88 & 80 & 32  \\
\ldots & \multicolumn{16}{l}{\ldots}
\end{array}
$$
\caption{Initial terms of sequence $\hat{t}(n)$ arranged in triangular form.}
\label{TabTPtri2}
\end{table}

The initial terms of the $\hat{t}(n)$ sequence are shown in
triangular array form in Table \ref{TabTPtri2}.
The initial $2^k-1$ terms of the $k$th row 
are the initial $2^k-1$ terms of the next row, so the rows
converge to a sequence 
$$
4,8,12,12,16,28,32,20,16,28,36,\ldots
$$
that we will denote by 
$F(n)$, $n \ge 0$ (this is entry A147646 in \cite{OEIS}).
It follows from the definition that $F(n)$ is defined
by the recurrence $F(0)=4$, $F(1)=8$, $F(2)=F(3)=12$, and for $n \ge 4$,
\beql{EqTP5}
F(2^k+i)=
\begin{cases}
2 F(i) + F(i+1),  &\text{if $0 \le i \le 2^k-3$}; \\
2 F(i) + F(i+1) -4,  &\text{if $i = 2^k-2$}; \\
2^{k+2}+4,  &\text{if $i = 2^k-1$}. \\
\end{cases}
\eeq
Also, for $0 \le i \le 2^k-2$,
\beql{EqTP6}
\hat{t}(2^k+i) = F(i) \,.
\eeq
Comparing \eqn{EqTP5} with Theorem \ref{th4},
we see (taking $\al=\be=4$, $\ga=1$, $\de=2$ in that theorem)
that $F(n)$ has generating function
\beql{EqTP7}
4(1+x) \prod_{k \ge 1} (1 + x^{2^k - 1} + 2 x^{2^k})
= 
2 \prod_{k \ge 0} (1 + x^{2^k - 1} + 2 x^{2^k}) \,.
\eeq
It follows from Theorem \ref{th6} that 
\beql{EqTP8}
F(i) = 
2 \sum_{m \ge 0} {2}^{\wt(i+m)-m} \binom{\wt(i+m)}{m} \,.
\eeq
Then \eqn{EqTP6} and \eqn{EqTP8} imply \eqn{EqHVE6}.~~~$\bsq$

%%%%%%%%%%%%%%%%%%%%%%%%%%%%%%%%%%%%%%
% Section: 
%%%%%%%%%%%%%%%%%%%%%%%%%%%%%%%%%%%%%%

\section{Cellular automata defined on graphs}\label{SecCAG}

Cellular automata defined on graphs were introduced 
by von Neumann and Ulam in 1949 \cite{Ulam50}, and
the toothpick structure (\S\ref{SecTP}),
the corner toothpick structure (\S\ref{SecCS}),
the leftist structure (\S\ref{SecLT}),
the Ulam-Warburton cellular automaton (\S\ref{SecUW}), etc.,
can all be described in this language. 
Let $G$ be an infinite directed (or undirected)
graph with finite in-degree and out-degree (or degree)
at each node. Nodes are either in the $\OFF$ state
or the $\ON$ state, and we just need to give a rule for 
deciding when the nodes change state.

For the toothpick structure, the nodes of $G$ are 
the vertices $(x,y) \in \ZZ \times \ZZ$ of
the square grid. There are two kinds of nodes:
``even'' nodes, with $x+y$ even, which have edges
directed to nodes $(x,y \pm 1)$,
and ``odd'' nodes, with $x+y$ odd, which have edges
directed to nodes $(x\pm 1,y)$.
Initially all nodes are $\OFF$, at stage $1$ we turn
$\ON$ node $(0,0)$, and thereafter a node turns $\ON$ 
if it has an incoming edge from exactly one $\ON$ node.
This is easily seen to be equivalent to the toothpick structure,
with the even (resp. odd) nodes representing 
the midpoints of vertical (resp. horizontal)
toothpicks.
The induced subgraph of $G$ joining the $\ON$ nodes
is a directed tree (and remains a tree if the arrows
on the edges are removed).

For the Ulam-Warburton cellular automaton, of course, 
$G$ is the undirected graph with vertices $(x,y) \in \ZZ \times \ZZ$,
with each node connected to its four neighbors
(or, in the case of the cellular automaton analyzed in
\S\ref{Sec8N}, its eight neighbors).
As already mentioned in \S\ref{SecUW}, the induced
subgraph joining the $\ON$ nodes is also a tree.

For the natural generalization of the Ulam-Warburton cellular automaton
to higher dimensions,
with $G= {\ZZ}^d$, $g \ge 1$,
and each node adjacent to its $2d$ neighbors,
there is an analog of Theorem \ref{th3}.
The number of cells that turn from $\OFF$ to $ON$ at
stage $n$ is given by $u(0)=0$, $u(1)=1$, and, for $n \ge 2$,
\beql{EqUWd}
u(n) = 2d \, (2d-1)^{\wt(n-1)-1}
\eeq
(Stanley and Chapman \cite{Stan94}; see also entries A151779, A151781 in \cite{OEIS}).
On the other hand, for the face-centered cubic lattice 
graph, where each node has $12$ neighbors, there is no obvious
formula (see A151776, A151777).

The number of possibilities is endless: see \cite{toothlist} and \cite{NKS} for
further examples.

\begin{figure}[htbp]
\centerline{\includegraphics[width=3.5in]{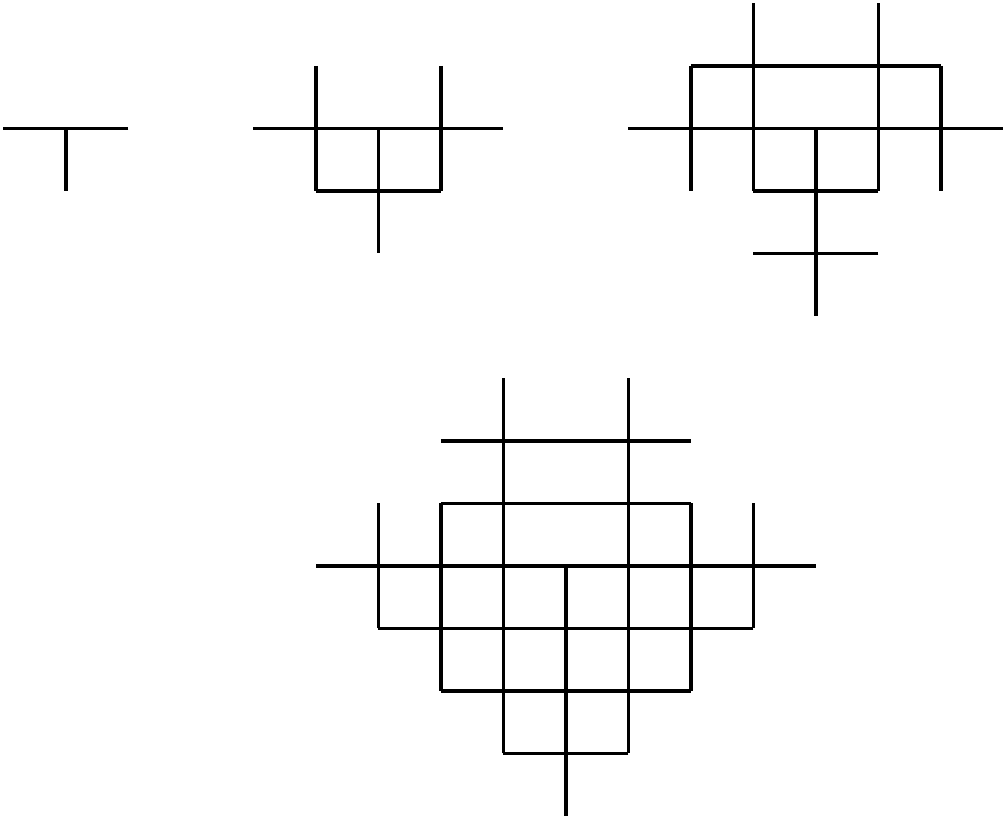}}
\caption{Stages $1$ through $4$ of the
evolution of the T-toothpick structure.}
\label{FigTT1}
\end{figure}

%%%%%%%%%%%%%%%%%%%%%%%%%%%%%%%%%%%%%%
% Section: 
%%%%%%%%%%%%%%%%%%%%%%%%%%%%%%%%%%%%%%

\section{T-shaped toothpicks}\label{SecTTP}

One may consider ``toothpicks'' of many other shapes.
Here we consider just one example: we define a ``T-toothpick''
to consist of three line segments of length $1$, forming a ``T''.
Using an obvious terminology, we will speak of the ``crossbar''
(which has length $2$)
and the ``stem'' (of length $1$) of the T, 
and refer to the midpoint of the crossbar
as the ``midpoint.'' 
A T-toothpick has three endpoints, and an endpoint is said to be 
``exposed'' if it is not the midpoint or endpoint of any other T-toothpick.

We start at stage $1$ with a single T-toothpick in the plane,
with its stem vertical and pointing downwards.
Thereafter, at stage $n$ we place a T-toothpick
at every exposed endpoint, with the midpoint of the new
T-toothpick touching the exposed endpoint and with its stem pointing
away from the existing T-toothpick.
Figure \ref{FigTT1} shows the first four stages of the evolution.
Let $\tau (n)$ denote the number of T-toothpicks added to
the structure at the $n$th stage (this is A160173).

\begin{theorem}\label{thTTP}
We have $\tau(0)=0$, $\tau(1)=1$, $\tau(2)=3$, and, for $n \ge 3$,
\beql{EqTTP1}
\tau(n) = \frac{2}{3}\{3^{\wt(n-1)}+3^{\wt(n-2)}\}+1 \,.
\eeq
\end{theorem}

\noindent{\bf Proof.} 
This is easily established from Theorem \ref{th3} by observing that 
the structures in the four quadrants defined by the the initial T are
essentially equivalent to copies of the Ulam-Warburton structure.
The copies in the first and second quadrants are one stage behind
those in the third and fourth quadrants.~~~$\bsq$

The analogous structure using Y-shaped toothpicks has resisted
our attempts to analyze it---see A160120 in \cite{OEIS}.

%%%%%%%%%%%%%%%%%%%%%%%%%%%%%%%%%%%%%%
% Section: 
%%%%%%%%%%%%%%%%%%%%%%%%%%%%%%%%%%%%%%

\section{The ``Maltese cross'' or Holladay-Ulam structure}\label{SecMC}

\begin{figure}[htbp]
\centerline{\includegraphics[width=2.5in]{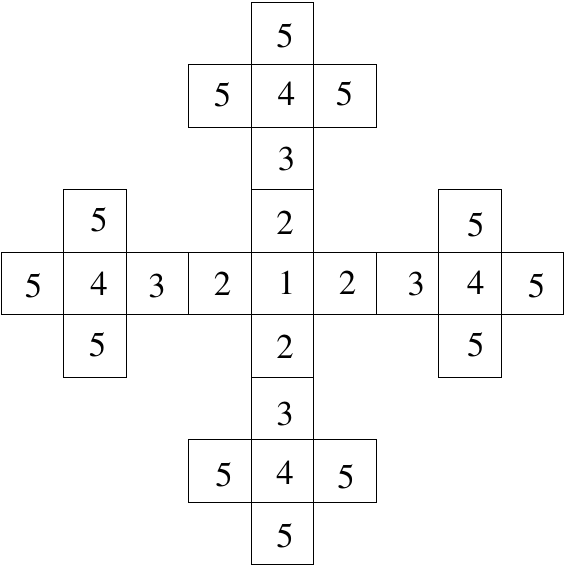}}
\caption{Stages $1$ through $5$ of the
evolution of the Maltese cross structure.}
\label{FigMalt1}
\end{figure}

On page 217 of \cite{Ulam62}, Ulam discusses another structure
that he and J. C. Holladay had studied.
To construct this, one first builds the
infinite Ulam-Warburton structure described in \S\ref{SecUW},
and then replaces each $\ON$ cell by a Maltese cross 
consisting of a central $\ON$ cell surrounded by 
four other $\ON$ cells.
Now label the cells of the new infinite structure,
starting by labeling the central square $1$, 
then the four adjacent cells $2$, and so on, 
always moving outwards from the center.
Figure \ref{FigMalt1} shows the cells with labels $1$ through $5$.

Holladay and Ulam then give a set of rules for a 
cellular automaton that will build up this structure,
starting with the cell labeled $1$. 
We found their rules (stated on pages 216 and 222)  somewhat hard to understand,
so it may be helpful to the reader if we give our version 
of them here.

For a candidate cell to be turned $\ON$ it is necessary (though
not sufficient) that it
share exactly one edge with an $\ON$ cell of the previous generation.
The two vertices of the candidate cell that touch that edge we will
call its {\em inner} vertices, and the other two vertices 
we call its {\em outer} vertices.

The rules are as follows. 
Cells are in one of three states,
$\OFF$, $\ON$ or $\DEAD$. 
Once a cell is $\ON$ or $\DEAD$ it remains in that state.
Initially all cells are $\OFF$, and at stage $1$ 
a single cell is turned $\ON$.

At stage $n$, consider $\OFF$ cells $X$ (the ``candidates'')
that share at least one edge with an $\ON$ cell of the previous generation.
If a candidate $X$ is edge-adjacent to two $\ON$ cells, $X$ is declared $\DEAD$.
If two candidates $X$ and $X'$ share an
outer vertex, both $X$ and $X'$ are declared $\DEAD$.
If a candidate $X$ is edge-adjacent to a $\DEAD$ cell, 
then $X$ is declared $\DEAD$,
except that when $n \equiv 2$ (mod $3$),
this rule obtains only if 
$X$ is edge-adjacent to a cell that was declared $\DEAD$
at the previous generation.
If a candidate is not excluded by these conditions,
it is turned $\ON$.

\begin{figure}[htbp]
\centerline{\includegraphics[width=2.5in]{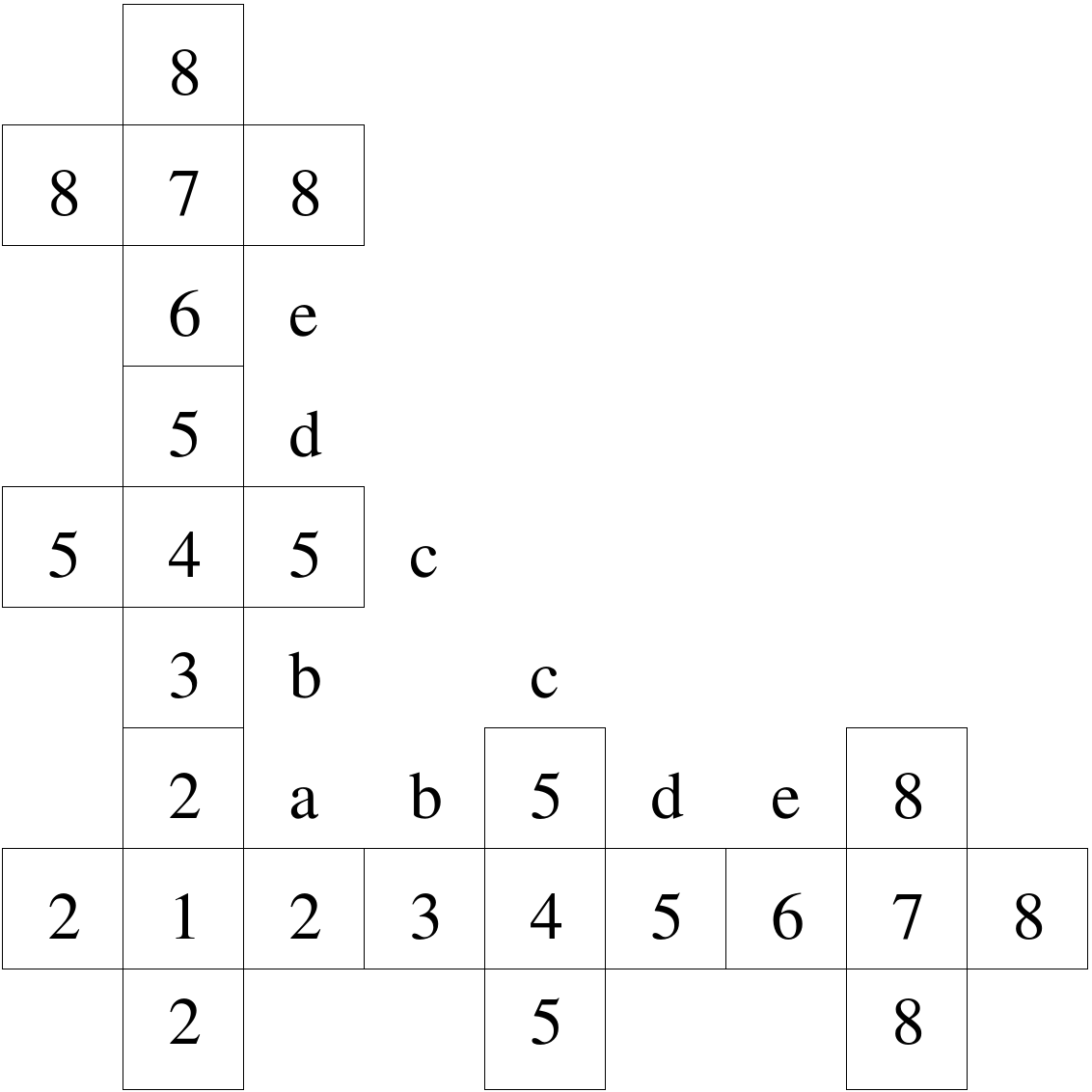}}
\caption{Annotated version of stages $1$ through $8$ of the
Maltese cross structure in the first quadrant.}
\label{FigMalt2}
\end{figure}

Figure \ref{FigMalt2} is an annotated picture showing the
first eight stages of the evolution of
this structure in the first quadrant, with letters identifying
the first few $\DEAD$ cells.
The letters 
`a' and `d' indicate candidate cells that are 
$\DEAD$ because they are edge-adjacent to
two $\ON$ cells, and
`b' and `e' indicate candidate cells that are
$\DEAD$ because they are edge-adjacent to
a $\DEAD$ cell and $n$ is not congruent to $2$ (mod $3$).
The `5' cells not on the main axes are $\ON$ because $5 \equiv 2$ (mod $3$).
The two `c' cells are $\DEAD$ because they would share
a common outer vertex; and so on.

Of course, the  original definition of this structure makes it easy
to give a formula for the number of cells, $m(n)$, say,
that are added at the $n$th stage. From \eqn{EqUW2} we have
$m(0)=0$, $m(1)=1$, $m(2)=4$, and, for $t \ge 1$,
\beql{EqMC1}
m(3t)=m(3t+1)=4 \mydot 3^{\wt(t)-1}, ~ m(3t+2)=4 \mydot 3^{\wt(t)}.
\eeq
(This is entry A151906 in \cite{OEIS}.)

We conclude this section by listing some related 
cellular automata studied by Ulam and his
colleagues Holladay and Schrandt in \cite{SU67}
and \cite{Ulam62}. 
Since we have not even found recurrences for them,
we will give no details.

Reference \cite{SU67} mentions a cellular automaton
which is intermediate in complexity between
the Ulam-Warburton structure and the Maltese cross
structure discussed above. This Schrandt-Ulam 
cellular automaton is described in
entries A170896 and A170897 in \cite{OEIS};
another version is given in A151895, A151896. 

Analogous structures based on triangles or hexagons
rather than squares are described in Examples $3$ through $6$
of \cite{Ulam62} (see also Wolfram \cite[p.~371]{NKS}), 
and in entries A151723, A151724, A161644, A161645
of \cite{OEIS}. We invite the reader to find recurrences
or generating functions for any of them.

%%%%%%%%%%%%%%%%%%%%%%%%%%%%%%%%%%%%%%
% Section: 
%%%%%%%%%%%%%%%%%%%%%%%%%%%%%%%%%%%%%%

\section{Rule 942}\label{Sec942}

On page 928 of \cite{NKS}, Wolfram considers (among other examples)
what happens if the rule for the Ulam-Warburton cellular automaton
is modified so that a cell turns $\ON$ if and only if either exactly one
or all four of its four neighbors is $\ON$
(this is ``Rule 942'' in the notation of \cite{PaWo85}, \cite{NKS}).

Let $w(n)$ denote the number of
cells that are changed from $\OFF$ to $\ON$ at stage $n$.
Since the four-neighbor part of the rule is invoked only
after an $\OFF$ cell is completely surrounded by $\ON$ cells, 
$w(n) \ge u(n)$ for all $n$.
In fact, $w'(n) := w(n)-u(n)$ is always a multiple of $4$,
and $w(n) = u(n)$ except when $n \equiv 1~(\mbox{mod~}n)$.
Table \ref{Tab942} shows the initial
values of $w(n)$, $u(n)$, $w'(n)$ and
$\de(n) := \frac{1}{4}(w(4n+1)-u(4n+1))$
(cf. A169648 and A169689 in \cite{OEIS}).

\begin{table}[htb]
$$
\begin{array}{|c|rrrrrrrrrrrrrrrr|}
\hline
n   & 0 & 1 & 2 & 3 & 4 & 5 & 6 & 7 & 8 & 9 & 10 & 11 & 12 & 12 & 14 & 15 \\
w(n) & 0 & 1 & 4 & 4 & 12 & 8 & 12 & 12 & 36 & 28 & 12 & 12 & 36 & 28 & 36 & 36 \\
u(n) & 0 & 1 & 4 & 4 & 12 & 4 & 12 & 12 & 36 & 4 & 12 & 12 & 36 & 12 & 36 & 36 \\
w'(n) & 0 & 0 & 0 & 0 & 0 & 4 & 0 & 0 & 0 & 24 & 0 & 0 & 0 & 16 & 0 & 0 \\
\de(n) & 0 & 1 & 6 & 4 & 24 & 4 & 20 & 12 & 84 & 4 & 20 & 12 & 76 & 12 & 60 & 36 \\
\hline
\end{array}
$$
\caption{The sequences $w(n)$, $u(n)$, $w'(n0 := w(n)-u(n)$, $\de(n)$
for $0 \le n \le 15$.}
\label{Tab942}
\end{table}

There is a simple explicit formula for $\de(n)$ and hence, via
\eqn{EqUW2}, for $w(n)$.

\begin{theorem}\label{th942}
We have $\de(0)=0$, $\de(1)=1$, $\de(2)=6$. For $n \ge 3$,
let $n = 2^k+j$ with $1 \le j \le 2^k$, where $j=2^m(2l+1)$ (say).
Then
\beql{Eq9421}
\de(n) = 4(3^{m+1}-2^{m+1})3^{\wt(l)} \,,
\eeq
except that if $j=2^k$ then
\beql{Eq9422}
\de(n) = 4 \mydot 3^{k+1}-3 \mydot 2^{k+1} \,.
\eeq
\end{theorem}
We omit the proof, which is similar to those
of Theorems \ref{th1} and \ref{th2}.

%%%%%%%%%%%%%%%%%%%%%%%%%%%%%%%%%%%%%%
% Section: 
%%%%%%%%%%%%%%%%%%%%%%%%%%%%%%%%%%%%%%

\section{Square grid with eight neighbors}\label{Sec8N}
Our final example is also based on the Ulam-Warburton cellular
automaton, except that now we take the neighbors of a cell
to consist of the eight cells surrounding it. 
(This is the Moore neighborhood of the cell, in
the notation of \cite{Mitch96}.)
Otherwise the rules are the same as in \S\ref{SecUW}:
a cell turns $\ON$ if exactly one of its eight neighbors is $\ON$.

\begin{figure}[htbp]
%\centerline{\includegraphics[width=5in]{A153006.pdf}}
\centerline{\includegraphics[width=4.0in]{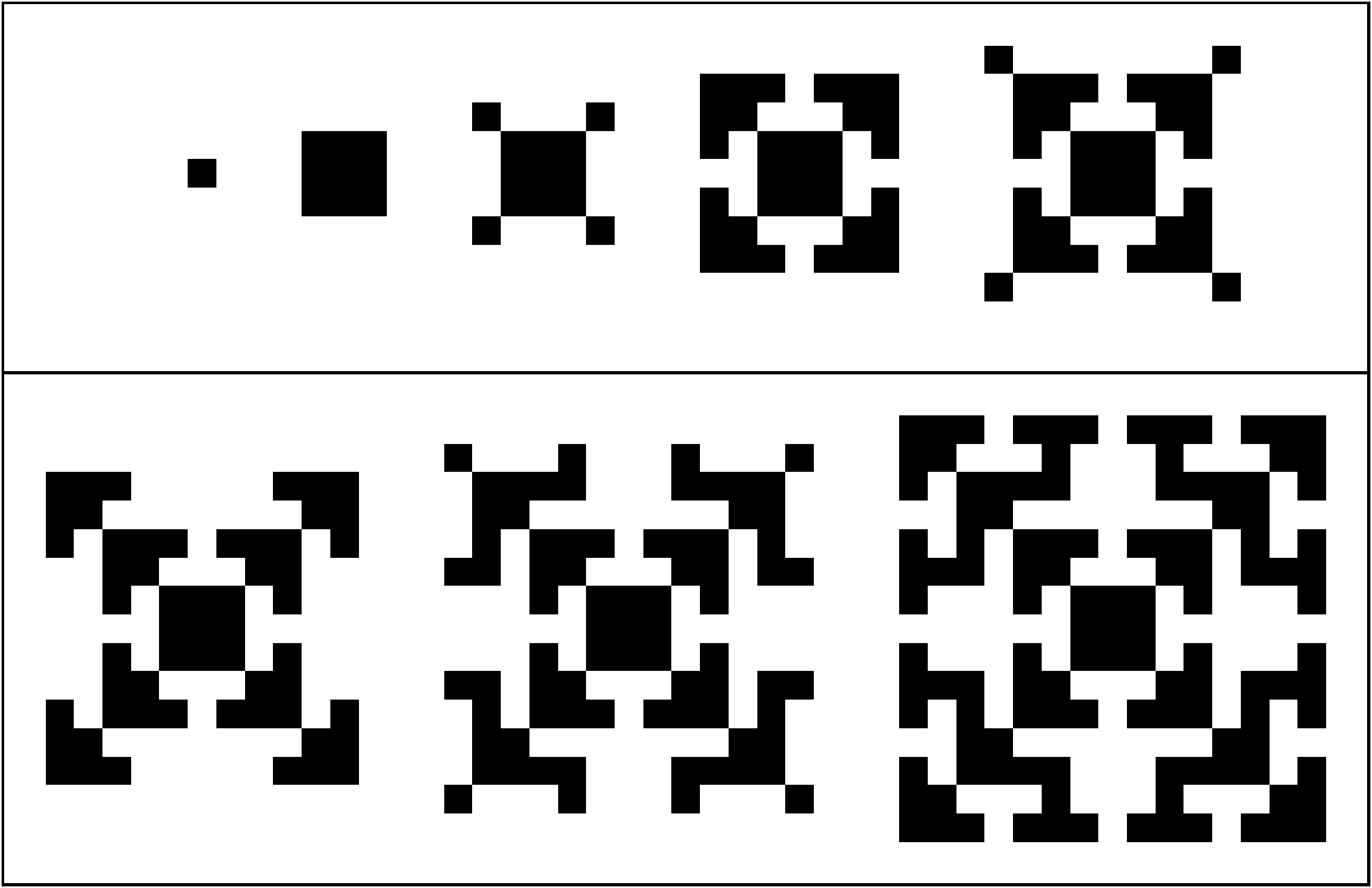}}
\caption{Stages $1$ through $8$ of the
evolution of the eight-neighbor structure.
The numbers of $\ON$ cells in the successive stages, $V(1), \ldots,
V(8)$, are $1, 9, 13, 33, 37, 57, 77, 121$.}
\label{Fig8NUnion}
\end{figure}

Let $v(n)$ ($n\ge 0$) denote the number of
cells that are changed from $\OFF$ to $\ON$ at the $n$th stage of
the evolution,
and let $V(n) := \sum_{i = 0}^{n} v(i)$ be the total number
of $\ON$ cells after $n$ stages.
The initial values of $v(n)$ and $V(n)$ are shown in Table \ref{Tab8N1}
below.
These sequences are respectively entries A151726 and A151725 in \cite{OEIS}.
Figure \ref{Fig8NUnion}
shows stages $1$ through $8$ of the
evolution of the this structure.

Since each cell now has two kinds of neighbors, it is
perhaps not surprising that this problem is more difficult
to analyze than the Ulam-Warburton structure.
In order to understand the growth, it is convenient to
define two versions of ``corner sequences,''
analogous to that introduced in \S\ref{SecCS}. 

So that we can refer to individual cells, we will
label each square cell by the grid point at its upper left corner.
That is, we define cell $(i,j)$ to consist of the
square $\{(x,y) \in \RR \times \RR \mid i \le x \le i+1,
~j-1 \le y \le j \}$.

\begin{figure}[htbp]
%\centerline{\includegraphics[width=5in]{A153006.pdf}}
\centerline{\includegraphics[width=4.5in]{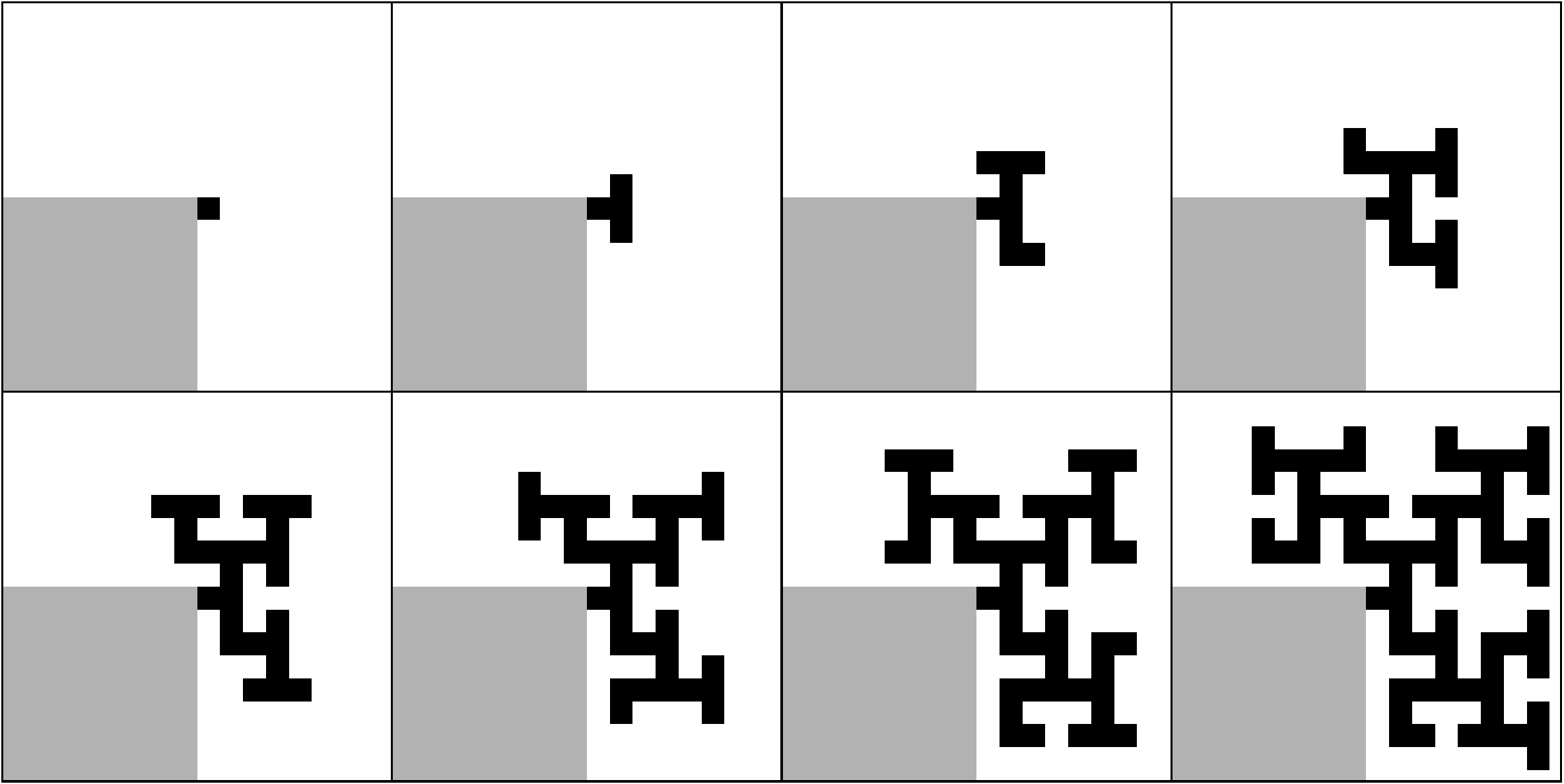}}
\caption{Stages $1$ through $8$ of the
evolution of the first corner sequence $v_1(n)$.}
\label{Fig8NUnionV1}
\end{figure}

For the first corner sequence, we exclude the third quadrant
of the plane, and at stage $1$ we turn $\ON$ the cell immediately
to the right of that quadrant (see Fig.~\ref{Fig8NUnionV1}).
More precisely, at stage $1$, we turn $\ON$ the cell $(0,0)$,
and thereafter extend the structure using the eight-neighbor rule,
with the proviso that after the first stage, no
$\ON$ cell may be adjacent to any
of the third-quarter cells---meaning the cells $(i,j) \in \ZZ \times \ZZ$
with $i \le -1, j \le 0$.

\begin{figure}[htbp]
%\centerline{\includegraphics[width=5in]{A153006.pdf}}
\centerline{\includegraphics[width=4.5in]{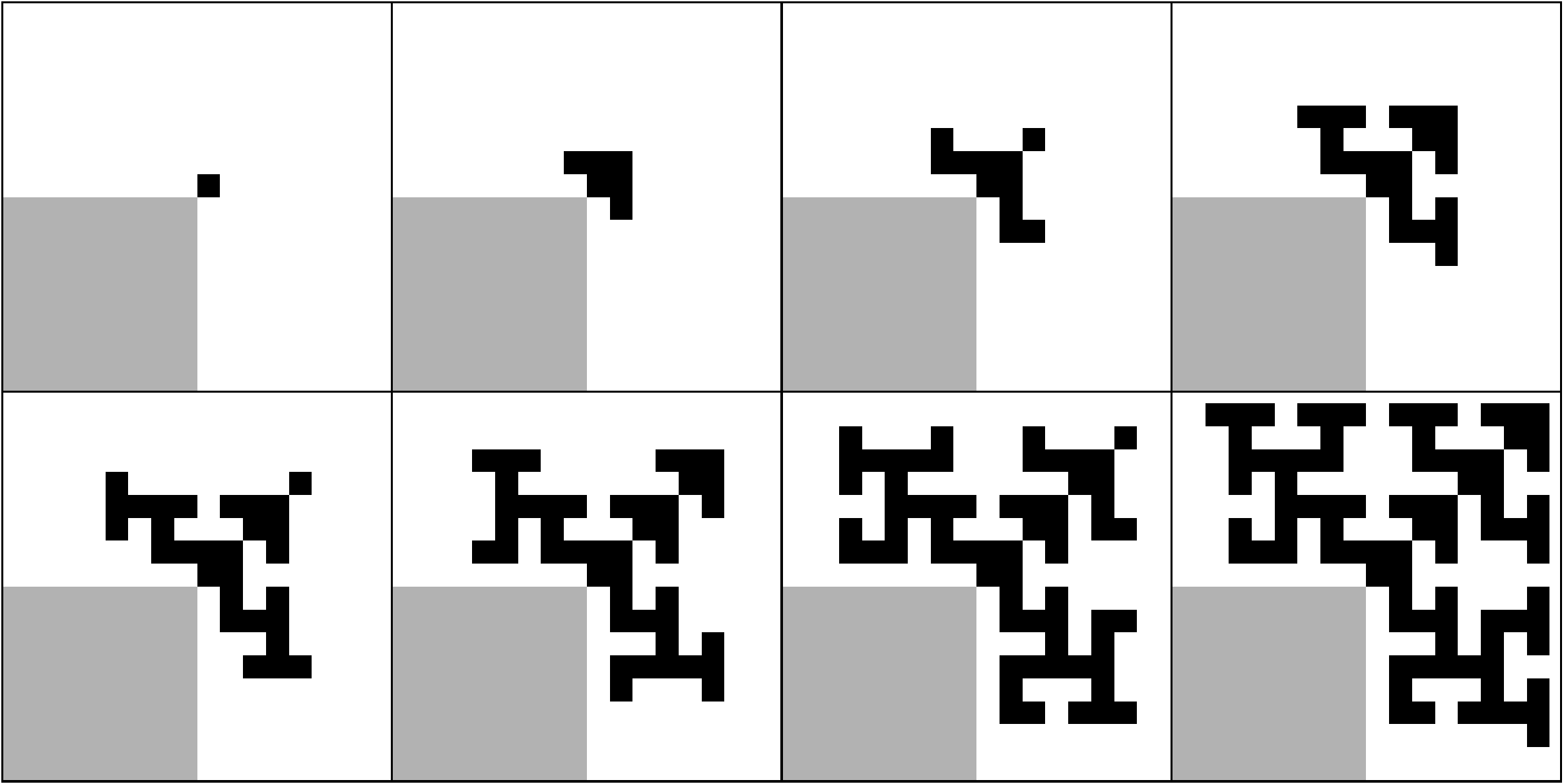}}
\caption{Stages $1$ through $8$ of the
evolution of the second corner sequence $v_2(n)$.}
\label{Fig8NUnionV2}
\end{figure}

The second corner sequence is similar to the first, except that at 
stage $1$ we turn $\ON$ the cell $(0,1)$, 
just up and to the right of the excluded quadrant
(Fig.~\ref{Fig8NUnionV2}).

Let $v_1(n)$ (resp. $v_2(n)$) denote the number of
cells that are changed from $\OFF$ to $\ON$ at the $n$th stage of
the evolution of the first (resp. second) corner sequence.
The initial values of $v_1(n)$ and $v_2(n)$ are also
shown in Table \ref{Tab8N1}.
These sequences are respectively entries A151747 and A151728 in \cite{OEIS}.
Figures \ref{Fig8NUnionV1} and \ref{Fig8NUnionV2}
shows stages $1$ through $5$ of the
evolution of the two corner sequences.

\begin{table}[htb]
$$
\begin{array}{|c|rrrrrrrrrr|}
\hline
n   & 0 & 1 & 2 & 3 & 4 & 5 & 6 & 7 & 8 & 9 \\
v(n) & 0 & 1 & 8 & 4 & 20 & 4 & 20 & 20 & 44 & 4 \\
V(n) & 0 & 1 & 9 & 13 & 33 & 37 & 57 & 77 & 121 & 125 \\
v_1(n) & 0 & 1 & 3 & 5 & 8 & 9 & 11 & 17 & 21 & 15 \\
v_2(n) & 0 & 1 & 5 & 5 & 11 & 7 & 15 & 19 & 23 & 7 \\
\hline
n & 10 & 11 & 12 & 13 & 14 & 15 & 16 & 17 & 18 & 19 \\
v(n) & 20 & 20 & 44 & 28 & 60 & 76 & 92 &  4 & 20 & 20 \\
V(n) & 145 & 165 & 209 & 237 & 297 & 373 & 465 & 469 & 489 & 509 \\
v_1(n) & 11 & 18 & 25 & 29 & 39 & 54 & 53 & 27 & 11 & 18 \\
v_2(n) & 15 & 21 & 29 & 29 & 49 & 59 & 47 & 7 & 15 & 21 \\
\hline
n & 20 & 21 & 22 & 23 & 24 & 25 & 26 & 27 & 28 & 29 \\
v(n) & 44 & 28 & 60 & 76 & 92 & 28 & 60 & 84 & 116 & 116 \\ 
V(n) & 553 & 581 & 641 & 717 & 809 & 837 & 897 & 981 & 1097 & 1213 \\
v_1(n) & 25 & 29 & 39 & 55 & 57 & 41 & 40 & 61 & 79 & 97  \\
v_2(n) & 29 & 29 & 49 & 61 & 53 & 29 & 51 & 71 & 87 & 107 \\
\hline
\end{array}
$$
\caption{The $8$-neighbor sequences $v(n)$ and $V(n)$,
and the two ``corner'' sequences $v_1(n)$, $v_2(n)$,
for $0 \le n \le 23$.}
\label{Tab8N1}
\end{table}

The following theorem gives recurrences for all three of
these sequences. 

\begin{theorem}\label{th8N}
The eight-neighbor sequences $v_1(n)$, $v_2(n)$ and $v(n)$ satisfy
the following recurrences:\\

$v_1(0)=0$, $v_1(1)=1$, $v_1(2)=3$, $v_1(3)=5$,
and, for $k \ge 2$,
\beql{Eq8NV1}
v_1(2^k+i)=
\begin{cases}
(3k+1)\, 2^{k-2}+1,  &\text{if $i=0$}; \\
3\mydot 2^{k-1}+3,         &\text{if $i=1$}; \\
2\, v_1(i) + v_1(i+1),     &\text{if $i=2, \ldots,2^k-2$}; \\
2\, v_1(i) + v_1(i+1)-1,   &\text{if $i=2^k-1$}; \\
\end{cases}
\eeq

$v_2(0)=0$, $v_2(1)=1$,
and, for $k \ge 1$,
\beql{Eq8NV2}
v_2(2^k+i)=
\begin{cases}
3\mydot 2^k-1,         &\text{if $i=0$}; \\
v_2(i) + 2\, v_1(i+1),     &\text{if $i=1, \ldots,2^k-2$}; \\
v_2(i) + 2\, v_1(i+1)-2,   &\text{if $i=2^k-1$}; \\
\end{cases}
\eeq

$v(0)=0$, $v(1)=1$,
and, for $k \ge 1$,
\beql{Eq8NV3}
v(2^k+i)=
\begin{cases}
6\mydot 2^k-4,         &\text{if $i=0$}; \\
4\, v_2(i),     &\text{if $i=1, \ldots,2^k-1$}. \\
\end{cases}
\eeq
\end{theorem}

Again we omit the proof.  We have not found generating 
functions or explicit formulas for any of these sequences.

\section{Acknowledgments}
We thank Beno\^{i}t Jubin for telling us about
his investigations into the asymptotic behavior of $T(n)$
that were discussed in \S\ref{SecFrac}, and
Gary Adamson and Hagen von Eitzen for their contributions to \cite{OEIS}
which were mentioned in \S\ref{SecGF} and \S\ref{SecEX}.
We also thank 
Maximilian Hasler,
John Layman and 
Richard Mathar,
who have made many contributions to \cite{OEIS}
(especially new sequences or extensions of existing sequences)
related to the subject of this paper. 
Finally, we thank Laurinda Alcorn for locating a copy of \cite{SU67}.

\end{document}